\documentclass[12pt]{article}

\usepackage{wrapfig}
\usepackage{graphicx}
\usepackage[none]{hyphenat}
\usepackage[english]{babel}
\usepackage[utf8]{inputenc}
\usepackage[T1]{fontenc}
\usepackage{parskip}

\usepackage{authblk}
\usepackage{tikz-cd}

\usepackage{caption}
\usepackage{csquotes}
\usepackage{mathtools,amssymb,amsfonts,amsthm}

\usepackage[breaklinks,unicode=true]{hyperref}
\usepackage[capitalise]{cleveref} 

\theoremstyle{plain}
\newtheorem{theorem}{Theorem}
\newtheorem*{theorem*}{Theorem}

\newtheorem*{lemma*}{Lemma}
\newtheorem{proposition}[theorem]{Proposition}
\newtheorem*{proposition*}{Proposition}

\newtheorem*{corollary*}{Corollary}

\theoremstyle{definition}
\newtheorem{definition}{Definition}
\newtheorem*{definition*}{Definition}
\newtheorem{example}{Example}
\newtheorem*{example*}{Example}

\crefname{theorem}{Theorem}{Theorems}
\Crefname{theorem}{Theorem}{Theorems}
\crefname{lemma}{Lemma}{Lemmas}
\Crefname{lemma}{Lemma}{Lemmas}
\crefname{proposition}{Proposition}{Propositions}
\Crefname{Prop}{Proposition}{Propositions}
\crefname{corollary}{Corollary}{Corollaries}
\Crefname{corollary}{Corollary}{Corollaries}
\crefname{definition}{Definition}{Definitions}
\Crefname{definition}{Definition}{Definitions}
\crefname{example}{Example}{Examples}
\Crefname{example}{Example}{Examples}
\crefname{theorem}{Theorem}{Theorems}

\newtheorem*{exercise*}{Exercise}
\crefname{exercise}{exercise}{exercises}
\Crefname{exercise}{Exercise}{Exercises}  

\theoremstyle{remark}
\newtheorem{remarkx}{Remark}
\newtheorem*{remarkx*}{Remark}
\crefname{remark}{Remark}{Remarks}
\Crefname{remark}{Remark}{Remarks}
\newenvironment{remark}
  {\pushQED{\qed}\remarkx}
  {\popQED\endremarkx}
\newenvironment{remark*}
  {\pushQED{\qed}\remarkx*}
  {\popQED\endremarkx*}

\title{New notions of uniformity and homogeneity of Cosserat media}

\author[-]{V\'ictor Manuel Jim\'enez}
\author[+*]{Manuel de León}

\affil[1]{\href{mailto:victor.jimenez@mat.uned.es}{victor.jimenez@mat.uned.es}}

\affil[2]{\href{mailto:mdeleon@icmat.es}{mdeleon@icmat.es}}

\affil[-]{Universidad Nacional de Educación a Distancia (UNED), Departamento de Matem\'aticas Fundamentales.
Calle de Juan del Rosal 10, 28040, Madrid, Spain} 

\affil[*]{Instituto de Ciencias Matem\'aticas (CSIC-UAM-UC3M-UCM),
C\textbackslash Nicol\'as Cabrera, 13-15, Campus Cantoblanco, UAM
28049 Madrid, Spain} 

\affil[+]{Real Academia de Ciencias Exactas, Fisicas y Naturales, C/de Valverde
22, 28004 Madrid, Spain}

\date{\today}

\begin{document}

\sloppy

\maketitle

\begin{abstract}

In this paper, we study internal properties of a Cosserat media. In  fact, by using groupoids and smooth distributions, we obtain a three canonical equations. The \textit{non-holonomic material equation for Cosserat media} characterizes the uniformity of the material. The \textit{holonomic material equation for Cosserat media} permits us to study when a Cosserat material is a second-grade material. It is remarkable that these two equations also provide us a unique and maximal division of the Cosserat medium into uniform and second-grade parts, respectively. Finally, we present a proper definition of homogeneity of the Cosserat medium, without assuming uniformity. Thus, the \textit{homogeneity equation for Cosserat media} characterizes this notion of homogeneity.
\end{abstract}

\tableofcontents

\section{Introduction}

\indent{In Continuum Mechanics a body is represented by a three-dimensional manifold $\mathcal{B}$ which can be covered with just one chart. A configuration $\psi$ is an embedding of $B$ in $\mathbb{R}^3$, and it is usual to identify the body with one of its configurations $\psi_0$, which is called a reference configuration. A change of configuration $\psi \circ \psi_0^{-1}$ is said to be a deformation.
 
A relevant problem is the following: given a mechanical response as a function of the positions on the body and the $1-$jets of the local diffeomorphisms of the body, how to decide if the body is uniform, namely, all the points of the body are made of the same material. A second crucial question is about the homogeneity of the body, which expresses the absence of defects of the material body.
 
W. Noll developed a geometric theory to deal with the properties of uniform bodies in his thesis \cite{WNOLLTHE} (see also \cite{WNOLL,CTRUE}) as well as its homogeneity. In that theory, he studies the model of \textit{simple materials} in which the properties of the body are encoded in the constitutive law, a function depending on the gradient of the deformation. This allows us to introduce the concept of material isomorphisms and re-interpret the uniformity in terms of the existence of a parallelism of such material isomorphisms. In addition, the homogeneity is characterized by the integrability of such parallelisms. Moreover, the use of $G-$structures has redefined the formulation and facilited the derivation of specific results (see for example \cite{MELZA,MAREMDL3,MAREMDL2,MAREMDL4,EPSBOOK2}). In fact, the lack of integrability of the associated $G-$structure manifests the presence of inhomogeneities (such as dislocations). Thus, we may say that the theory of inhomogeneities of smoothly uniform simple materials is well established in terms of differential geometric structures. Nevertheless, in the absence of uniformity, these $G-$structures do not exist.
 
As a further step in the study of uniformity and homogeneity, one can use the notion of groupoid; indeed, the collection of all the possible material isomorphisms on a body $\mathcal{B}$ has the structure of subgroupoid of the Lie groupoid of all $1-$jets of local diffeomorphisms from $\mathcal{B}$ into itself, $\Pi^1\left( \mathcal{B} , \mathcal{B}\right)$. As a difference with the $G-$structures, the existence of this groupoid, called material groupoid, does not depend on the uniformity of the material. In fact, the material groupoid may always be constructed for any material body.

However, the material groupoid is not necessarily differentiable. Indeed, the property of differentiability of this groupoid works to characterize some material property, as for example the smooth uniformity. In a series of papers \cite{MGEOEPS,MD,CHARDIST,GENHOM}, we have developed a new theory of uniformity and homogeneity without the necessity of smoothness. Indeed, even if a subgroupoid of a Lie groupoid is not a Lie subgroupoid itself, one can generalize the construction of the asociated Lie algebroid, obtaining the so-called characteristic distributions. This result permits us to extend the concept of uniformity as well as consider homogeneity even if we are not in presence of uniformity.

However, there are many non-simple materials. In fact, materials like granular solids, rocks or bones cannot be modelled without extra kinematic variables \cite{GCAPR}. The theory of generalized media was introduced by Eugène and François Cosserat between 1905 and 1910. The Cosserats associated to each point of the body a family of vector (directors). In a more mathematical way, a Cosserat continuum can be described as a manifold of dimension $m$ and a family of $n$ vector fields on the manifold. Some of the developments of the theory can be found in Maugin \cite{GAMAU1,GAMAU3} or in Kröner \cite{EKRON}; we also remit to the work by Eringen \cite{ERIN}. A particular case of Cosserat material are the second-grade materials. According with the article \cite{MAREMDL}, $\mathcal{B}$ is a second-grade material, if all the material isomorphisms are natural prolongations to the frame bundle of the induced diffeomorphisms on the basis. In other words, all the material isomorphisms are $1-$jets $j_{X,Y}^{1} F \psi$, where $\psi$ is a local diffeomorphism on the body $\mathcal{B}$.
 
The geometrical structures which are necessary to develop a rigorous theory have been available for some time. Actually, the notion of director given by Cosserats is closely related with frame bundles. In 1950 C. Ehresmann (see \cite{CELP,CELC11,CELC12,CELC13}) formalized the notion of principal bundles and studied many frame bundles associated in a natural way to an arbitrary manifold: non-holonomic and holonomic frame bundles. Thus, we can intepret a Cosserat medium as a linear frame bundle $F\mathcal{B}$ of a manifold $\mathcal{B}$, the macromedium, which can be covered with just one chart (see \cite{EPSTEIN1998127}). Then, a configuration of $F\mathcal{B}$ is an embedding $\Psi : F\mathcal{B} \rightarrow F\mathbb{R}^{3}$ of principal bundles such that the induced Lie group morphism is the identity map. We fix a configuration $\Psi_{0}$, as the reference configuration, and a deformation is a change of configurations, $\chi = \Psi \circ \Psi_{0}^{-1}$.
 
The constitutive elastic law is now written as
$$ W \ = \ W\left( X , F \right),$$
where $X$ is a point of the macromedium and $F$ is the gradient of a deformation $\chi$ at a point $X$. Since $\chi$ is a morphism of principal bundles, $F$ depends only on the base points. This constitutive equation permits us to associate to each two points $X,Y$ of $\mathcal{B}$ the family of material isomorphisms from $X$ to $Y$ (which could be empty), i.e., the $1-$jets $G$ at $X$ of the local principal bundle isomorphisms from $X$ to $Y$ which satisfy that
$$ W\left( X , F \cdot G\right) \ = \ W\left( Y , F \right),$$
for all deformation gradients $F$ at $Y$.

The present paper is devoted to extend the construction of the characteristic distribution for Cosserat media. In particular, we have proved that, associated to any Cosserat media there are two well-defined characteristic distributions, the non-holonomic material distribution of second order and the holonomic material distribution of second order, respectively.

These distributions are generated by the left-invariant vector fields which are in the kernel of $TW$ (Eq. (\ref{14.22})) and the complete lift of vector fields on the macromedium $\mathcal{B}$ which are in the kernel of $TW$, respectively. These facts induce two distinct equations which permit us to construct these characteristic distribution,
\begin{itemize}
    \item \textbf{Non-holonomic material equation for Cosserat media} (\ref{materialeqcoss2345})
 \begin{equation*}
  \hspace{-1cm}  -\Theta^{i} \dfrac{\partial W }{\partial x^{i} } + \Theta^{l}_{i}\left[ y^{j}_{l}\dfrac{\partial  W}{\partial y^{j}_{i} } +  y^{j}_{l,k} \dfrac{\partial  W}{\partial y^{j}_{i,k} }\right] + \Theta^{l}_{,i} \left[ y^{j}_{,l}  \dfrac{\partial  W}{\partial y^{j}_{,i} } +  y^{j}_{m,l}\dfrac{\partial  W}{\partial y^{j}_{m,i} }\right] +   \Theta^{l}_{i,k} y^{j}_{l} \dfrac{\partial  W}{\partial y^{j}_{i,k} }= 0 
\end{equation*}

\item \textbf{Holonomic material equation for Cosserat media} (\ref{holonomicmaterialeqcoss2245})

\begin{small}
\begin{equation*}
\hspace{-1cm}-\Theta^{i} \dfrac{\partial W }{\partial x^{i} } + \dfrac{\partial  \Theta^{l}}{\partial x^{i} }\left[ y^{j}_{l}\dfrac{\partial  W}{\partial y^{j}_{i} }  +  y^{j}_{,l}  \dfrac{\partial  W}{\partial y^{j}_{,i} } + y^{j}_{l,k} \dfrac{\partial  W}{\partial y^{j}_{i,k} } +  y^{j}_{m,l}\dfrac{\partial  W}{\partial y^{j}_{m,i} }\right] +    \dfrac{\partial^{2} \Theta^{l} }{\partial x^{i}\partial x^{k} } y^{j}_{l} \dfrac{\partial  W}{\partial y^{j}_{i,k} }= 0 .
\end{equation*}
\end{small}

\end{itemize}

Furthermore, we have proved that the Cosserat media may be canonically divided into smoothly uniform \textit{parts}, and second-grade \textit{parts}. The spaces of solutions of these equations characterizes these two properties (Theorems \ref{uniformityth}, \ref{uniformityandsecondth}, \ref{uniformityandsecondth23}, and  \ref{uniformityandsecondth23456}).

Finally, we have studied the property of homogeneity. In particular, we defined, by first time, a notion of homogeneity which is valid for non-uniform Cosserat materials and generalizes the well-known notion of homogeneity. Roughly speaking, a Cosserat material will be homogeneous, when each smoothly uniform material part is homogeneous and all the uniform material submanifolds can be ``\textit{straightened at the same time}''. Thus, we found another differential equation,
\begin{itemize}
    \item \textbf{Homogeneity equation for Cosserat media} (\ref{2homogeneousmaterialeqcoss2345})
    \begin{equation*}
  \dfrac{\partial  W}{\partial x^{k} } + y^{j}_{l}\dfrac{\partial P^{l}_{i} }{\partial x^{k}} \dfrac{\partial  W}{\partial y^{j}_{i} } + \left(y^{j}_{l} \left[\dfrac{\partial^{2} P^{l}_{i} }{\partial x^{k}x^{m}} + \dfrac{\partial^{2} P^{l}_{i} }{\partial x^{k}y^{m}} \right]+ y^{j}_{l,m}\Theta^{l}_{i}\right) \dfrac{\partial  W}{\partial y^{j}_{i,m} } = 0,
\end{equation*}
\end{itemize}
which characterizes this intuitive notion of homogeneity.

The paper is structured as follows. Section \ref{Chardistbackg94} is devoted to present a brief introduction to frame bundles and groupoids, focusing on the properties and examples which are necessary for the development of the paper. Here, we recall the notion of characteristic distribution associated to an arbitrary subgroupoid of a Lie groupoid. Section \ref{CosseratSection234} is devoted to present Cosserat materials, next to the relevant definitions associated to this model. It is important to recall that these two sections are introductory. In fact, the novelty results start with section \ref{CosseratDistrib233}. Here, we present a detailed construction of the non-holonomic material distribution of second order, the holonomic material distribution of second order, the non-holonomic material equation for Cosserat media (\ref{materialeqcoss2345}), and the holonomic material equation for Cosserat media (\ref{holonomicmaterialeqcoss2245}). We prove here that the Cosserat media is uniquely divided into smoothly uniform \textit{parts}, and second-grade \textit{parts}. Furthermore, we the solutions of these equations permits us to calculate these divisions of the material body in a precise way.

Finally, section \ref{Homogeneitysection24} is devoted to define and study the homogeneity of a (non necessarily) uniform Cosserat body. In particular, the solutions of the homogeneity equation for Cosserat media (\ref{2homogeneousmaterialeqcoss2345}) permits us to characterize this new notion of homogeneity.

\section{On frame bundles and groupoids}\label{Chardistbackg94}

Let us start with the notion of \textit{principal bundle} (see \cite{KONOM}) and, as relevant cases, we will introduce the concept of \textit{frame bundle of a manifold}. This will permit us to work with \textit{Cosserat materials}.\\
\begin{definition}
\rm
Let $P$ be a manifold and $G$ be a Lie group which acts over $P$ by the right satisfying:

\begin{itemize}
\item[(i)] The action of $G$ is free, i.e., 
$$xg = x \Leftrightarrow g = e,$$
where $e \in G$ is the identity of $G$.
\item[(ii)] The canonical projection $\pi: P \rightarrow M = P/G$, where $P/G$ is the space of orbits, is a surjective submersion.
\item[(iii)] $P$ is locally trivial, i.e., for each point $x \in M$ there is a neighborhood $U$ of $x$ such that $P$ is locally a product $U \times G$. More precisely, there exists a diffeomorphism $\Phi : \pi^{-1}\left(U\right) \rightarrow U \times G$, such that $\Phi \left(u\right) = \left(\pi \left(u\right) , \phi\left(u\right)\right)$, where the map $\phi : \pi^{-1}\left(U \right) \rightarrow G$ satisfies that
$$\phi \left(ua \right) = \phi \left( u \right)a, \ \forall u \in U, \ \forall a \in G.$$
$\Phi$ is called a \textit{trivialization on $U$}.

\end{itemize}
\end{definition}
A principal bundle will be denoted by $P\left(M,G\right)$, or simply $\pi: P \rightarrow M$ if there is no ambiguity about to the structure group $G$. $P$ is called the \textit{total space}, $M$ is the \textit{base space}, $G$ is the \textit{structure group} and $\pi$ is the \textit{projection}. The closed submanifold $\pi^{-1}\left(x\right)$, $x \in M$ will be called the \textit{fibre over $x$}. For each point $u \in P$, we have $\pi^{-1}\left(x\right) =  uG$, where $\pi \left(u\right) = x$, and $u G$ will be called the \textit{fibre through $u$}. Every fibre is diffeomorphic to $G$, but this diffeomorphism depends on the choice of the trivialization.\\

\begin{definition}
\rm
Given $P\left(M,G\right)$ and $P'\left(M',G'\right)$ principal bundles, a principal bundle morphism from $P\left(M,G\right)$ to $P'\left(M',G'\right)$ consists of a differentiable map $\Phi: P \rightarrow P'$ and a Lie group homomorphism $\varphi : G \rightarrow G'$ such that
$$\Phi \left(ua\right) = \Phi\left(u\right) \varphi\left(a\right), \ \forall u \in P, \ \forall a \in G.$$
\end{definition}

\noindent Notice that, in this case, $\Phi$ maps fibres into fibres and it induces a differentiable map $\phi :M \rightarrow M'$ by the equality $ \phi\left(x\right) = \pi\left(\Phi\left(u\right)\right)$, where $u \in \pi^{-1}\left(x\right)$.\\
$P\left(M,G\right)$ is said to be a \textit{subbundle} of $P'\left(M',G'\right)$ in case that the maps characterizing the principal bundle morphism are embeddings. In such a case, we can identify $P$ with $\Phi\left(P\right)$, $G$ with $\varphi \left(G\right)$ and $M$ with $\phi\left(M\right)$.\\
Finally, a principal bundle morphism is called \textit{isomorphism} if it can be inverted by another principal bundle morphism.\\

\begin{example}
\rm
Given a manifold $M$ and $G$ a Lie group, we can consider $M \times G$ as a principal bundle over $M$ with projection $pr_{1} : M \times G \rightarrow M$ and structure group $G$. The action is given by,
$$\left(x,a\right)b = \left(x,ab\right),  \ \forall x \in M , \ \forall a,b \in G.$$
This principal bundle is called a \textit{trivial principal bundle}.
\end{example}
Now, we will introduce an important example of principal bundle, the \textit{frame bundle of a manifold}. In order to do that, we will start with the following definition.
\begin{definition}
\rm
Let $M$ be a manifold. A \textit{linear frame} at the point $x \in M$ is an ordered basis of $T_{x}M$.
\end{definition}

\begin{remark}
\rm
Alternatively, a linear frame at $x$ can be viewed as a linear isomorphism $z : \mathbb{R}^{n} \rightarrow T_{x}M$ identifying a basis on $T_{x}M$ as the image of the canonical basis of $\mathbb{R}^{n}$ by $z$.\\
There is a third way to interpret a linear frame by using the theory of jets. Indeed, a linear frame $z$ at $x \in M$ may be considered as the 1-jet $j^{1}_{0,x}\phi$ of a local diffeomorphism $\phi$ from an open neighbourhood of $0$ in $\mathbb{R}^{n}$ onto an open neighbourhood of $x$ in $M$ such that $\phi\left(0\right) = x$. So, $z= T_{0}\phi$.\\

\end{remark}
Thus, we denote by $FM$ the set of all linear frames at all the points of $M$. We can view $FM$ as a principal bundle over $M$ with the structure group $Gl\left(n , \mathbb{R}\right)$ and projection $\pi_{M} : FM \rightarrow M$ given by
$$ \pi_{M} \left(j^{1}_{0,x}\phi\right) = x, \ \forall j^{1}_{0,x}\phi \in FM.$$
This principal bundle is called the \textit{frame bundle on} $M$. Let $\left(x^{i}\right)$ be a local coordinate system on an open set $U \subseteq M$. Then we can introduce local coordinates $\left(x^{i}, x^{i}_{j}\right)$ over $FU \subseteq FM$ such that
\begin{equation}\label{27}
x^{i}_{j} \left( j_{0,x}^{1}\phi \right) =  \dfrac{\partial  \left(x^{i}\circ \phi\right)}{\partial x^{j}_{|0} }.
\end{equation}
If $\psi : N \rightarrow M$ is a local diffeomorphism, we denote by $F \psi : FN \rightarrow FM$ the local isomorphism induced from $\phi$, and defined by
$$F \psi \left( j_{0, x}^{1} \phi \right) = j_{0, \psi \left( x\right)}^{1}\left( \psi \circ \phi\right).$$
We will denote by $e_{1}$ the frame $j_{0,0}^{1}Id_{\mathbb{R}^{n}} \in F \mathbb{R}^{n}$, where $Id_{\mathbb{R}^{n}}$ is the identity map on $\mathbb{R}^{n}$. Let $ \Psi : F \mathbb{R}^{n}  \rightarrow FM$ be a local isomorphism of principal bundles such that its domain contains 
$e_{1}$, and whose induced isomorphism on Lie groups is the identity, i.e.,
$$\Psi \left( z \cdot g\right) = \Psi \left(z\right) \cdot g, \ \forall z \in Dom\left(\Psi\right) \subseteq F \mathbb{R}^{n}, \ \forall g \in Gl\left(n, \mathbb{R}^{n}\right).$$
We denote by $\psi : \mathbb{R}^{n} \rightarrow M$ the local diffeomorphism induced by $\Psi$, i.e.,
$$ \psi \circ \pi_{\mathbb{R}^{n}} = \pi_{M} \circ \Psi.$$
Notice that, $j_{e_{1}, \Psi \left( e_{1}\right)}^{1} \Psi$ can be identified with a linear frame at the point $\Psi \left(e_{1}\right)$ since $T_{e_{1}} \Psi  : T_{e_{1}} \left( F \mathbb{R}^{n} \right) \cong \mathbb{R}^{n+n^{2}} \rightarrow T_{\Psi \left( e_{1}\right)}FM$ is a linear isomorphism, and, therefore, the collection of all $1-$jets $j_{e_{1}, \Psi \left( e_{1}\right)}^{1} \Psi$, denoted by $\overline{F}^{2}M$, is a submanifold of $F\left(FM\right)$. Associated to this manifold, we may construct three canonical projections $\overline{\pi}^{2}_{1},\tilde{\pi}^{2}_{1} : \overline{F}^{2}M \rightarrow FM$ and $\overline{\pi}^{2} : \overline{F}^{2}M \rightarrow M$ given by:
\begin{itemize}
\item $\overline{\pi}^{2}_{1} \left(j_{e_{1} , \Psi \left( e_{1} \right)}^{1} \Psi\right) = \Psi \left(e_{1}\right)$
\item $\tilde{\pi}^{2}_{1} \left(j^{1}_{e_{1}, \Psi \left( e_{1} \right)} \Psi \right) = j_{0,z}^{1} \psi.$
\item $\overline{\pi}^{2} \left(j_{e_{1} , \Psi \left( e_{1} \right)}^{1} \Psi\right) = \psi \left(0\right)$
\end{itemize}
These projections are related as the following commutative diagram shows

\vspace{35pt}

\vspace{0.5cm}
\begin{picture}(375,50)(50,40)
\put(200,20){\makebox(0,0){$FM$}}
\put(250,25){$\pi_{M}$}               \put(220,20){\vector(1,0){70}}
\put(260,65){$\overline{\pi}^{2} $}
\put(220,100){\vector(1,-1){70}}
\put(310,20){\makebox(0,0){$M$}}
\put(165,50){$\overline{\pi}^{2}_{1}$}                  \put(200,95){\vector(0,-1){60}}
\put(320,50){$\pi_{M}$}                  \put(310,95){\vector(0,-1){60}}
\put(200,115){\makebox(0,0){$\overline{F}^{2}M$}}
\put(250,120){$\tilde{\pi}^{2}_{1}$}               \put(220,115){\vector(1,0){70}}
\put(310,115){\makebox(0,0){$FM$}}
\end{picture}

 \vspace{35pt}

A direct computation shows that $\overline{F}^{2}M$ is a principal bundle over $FM$ with canonical projection $\overline{\pi}^{2}_{1}$ and structure group,
$$ \overline{G}^{2}_{1}\left(n\right) := {\{ j_{e_{1} , e_{1}}^{1} \Psi \in \overline{F}^{2}\mathbb{R}^{n} / \ \Psi \left(e_{1}\right)= e_{1} \} = \overline{\pi}^{2}_{1}}^{-1} \left(e_{1}\right).$$
Notice that $\overline{G}^{2}_{1}\left(n\right)$ is a Lie subgroup of $Gl\left(n+n^{2} , \mathbb{R}\right)$ acting on $\overline{F}^{2}M$ by composition of jets.  We also have that $\overline{F}^{2}M$ is a principal bundle over $M$ with canonical projection $\overline{\pi}^{2}$ and structure group
$$ \overline{G}^{2}\left(n\right) := \{ j_{e_{1} , \Psi \left(e_{1}\right)}^{1} \Psi \in \overline{F}^{2}\mathbb{R}^{n} / \ \psi \left(0\right)= 0 \} = {\overline{\pi}^{2}}^{-1} \left(0\right),$$
which, again, acts on $\overline{F}^{2}M$ by composition of jets. The principal bundle $\overline{F}^{2} M$ will be called the \textit{non-holonomic frame bundle of second order} and its elements will be called \textit{non-holonomic frames of second order}. There are more principal bundles defined over the $1-$jets of local isomorphisms $j^{1}_{e_{1} , \Psi \left( e_{1} \right)} \Psi$ on $FM$ (\textit{holonomic} and \textit{semi-holonomic}). To know about the relations between them see \cite{MDELAM}.\\

By taking into account the coordinates defined on $FM$, given a local coordinate system $\left(x^{i}\right)$ on an open set $U \subseteq M$, we can introduce local coordinates $\left(x^{i}, x^{i}_{j}\right)$ over $FU \subseteq FM$ and, hence, we can also introduce local coordinates \linebreak$\left(\left(x^{i}, x^{i}_{j}\right),x^{i}_{,j},x^{i}_{,jk},x^{i}_{j,k}, x^{i}_{j,kl}\right)$ over $F\left(FU\right)$ such that
\begin{itemize}
\item $x^{i}_{,j} \left( j_{e_{1},Z}^{1}\Psi \right) =  \dfrac{\partial \left(x^{i}\circ \Psi\right)}{\partial x^{j}_{|e_{1}} }$
\item $x^{i}_{,jk} \left( j_{e_{1},Z}^{1}\Psi \right) =  \dfrac{\partial \left(x^{i}\circ \Psi\right)}{\partial {x^{j}_{k}}_{|e_{1}} }$
\item $x^{i}_{j,k} \left( j_{e_{1},Z}^{1}\Psi \right) =  \dfrac{\partial \left(x^{i}_{j}\circ \Psi\right)}{\partial {x^{k}}_{|e_{1}} }$
\item $x^{i}_{j,kl} \left( j_{e_{1},Z}^{1}\Psi \right) =  \dfrac{\partial \left(x^{i}_{j}\circ \Psi\right)}{\partial {x^{k}_{l}}_{|e_{1}} }$
\end{itemize}
Thus, restricting to $\overline{F}^{2}U$ we have that
\begin{itemize}
\item $x^{i}_{,jk} =  0$
\item $x^{i}_{j,kl} =x^{i}_{k}\delta_{l}^{j}$
\end{itemize}
Then, the induced coordinates on $\overline{F}^{2}U$ are given by
\begin{equation}\label{138.nonh}
\left(\left(x^{i}, x^{i}_{j}\right),x^{i}_{,j},x^{i}_{j,k}\right)
\end{equation}
in such a way that
\begin{itemize}
\item $\pi_{M}\left(x^{i}, x^{i}_{j}\right) = x^{i}$
\item $\pi_{FM}\left(\left(x^{i}, x^{i}_{j}\right),x^{i}_{,j},x^{i}_{,jk},x^{i}_{j,k}, x^{i}_{j,kl}\right) = \left(x^{i}, x^{i}_{j}\right)$
\item $\overline{\pi}^{2}_{1} \left(\left(x^{i}, x^{i}_{j}\right),x^{i}_{,j},x^{i}_{j,k}\right) = \left(x^{i} , x^{i}_{j}\right)$
\item $\overline{\pi}^{2} \left(\left(x^{i}, x^{i}_{j}\right),x^{i}_{,j},x^{i}_{j,k}\right) = x^{i}$
\item $\tilde{\pi}^{2}_{1} \left(\left(x^{i}, x^{i}_{j}\right),x^{i}_{,j},x^{i}_{j,k}\right) = \left(x^{i}, x^{i}_{,j}\right)$
\end{itemize}

We will give here a very brief introduction on \textit{(Lie) groupoids} and the relation with \textit{(smooth) distributions} which is crucial to understand the results proved in this paper. For a detailed study we refer to \cite{VMMDME} (see also \cite{CHARDIST,MD}). For groupoids we recommend \cite{KMG}.

\begin{definition}
\rm
Let $ M$ be a set. A \textit{groupoid} over $M$ is given by a set $\Gamma$ equipped with the maps $\alpha,\beta : \Gamma \rightarrow M$ (\textit{source map} and \textit{target map} respectively), $\epsilon: M \rightarrow \Gamma$ (\textit{section of identities}), $i: \Gamma \rightarrow \Gamma$ (\textit{inversion map}) and $\cdot : \Gamma_{\left(2\right)} \rightarrow \Gamma$ (\textit{composition law}). Here, $\Gamma_{(k)}$ denotes the $k$-tuplas $ \left(g_{1}, \hdots , g_{k}\right) \in \Gamma \times \stackrel{k)}{\ldots} \times \Gamma $ such that $\alpha\left(g_{i}\right)=\beta\left(g_{i+1}\right)$ for $i=1, \hdots , k -1$. The following properties are satisfied:\\
\begin{itemize}
\item[(1)] $\alpha$ and $\beta$ are surjective and for each $\left(g,h\right) \in \Gamma_{\left(2\right)}$,
$$ \alpha\left(g \cdot h \right)= \alpha\left(h\right), \ \ \ \beta\left(g \cdot h \right) = \beta\left(g\right).$$
\item[(2)] Associativity of the composition law, i.e.,
$$ g \cdot \left(h \cdot k\right) = \left(g \cdot h \right) \cdot k, \ \forall \left(g,h,k\right) \in \Gamma_{\left(3\right)}.$$
\item[(3)] For all $ g \in \Gamma$,
$$ g \cdot \epsilon \left( \alpha\left(g\right)\right) = g = \epsilon \left(\beta \left(g\right)\right)\cdot g .$$
In particular,
$$ \alpha \circ  \epsilon \circ \alpha = \alpha , \ \ \ \beta \circ \epsilon \circ \beta = \beta.$$
\item[(4)] For each $g \in \Gamma$,
$$i\left(g\right) \cdot g = \epsilon \left(\alpha\left(g\right)\right) , \ \ \ g \cdot i\left(g\right) = \epsilon \left(\beta\left(g\right)\right).$$
Then,
$$ \alpha \circ i = \beta , \ \ \ \beta \circ i = \alpha.$$
\end{itemize}
These maps ($\alpha$, $\beta$, $\epsilon$, $i$, and $\cdot$) will be called the \textit{structure maps}. We will denote this groupoid by $ \Gamma \rightrightarrows M$.
\end{definition}
Observe that, since $\alpha$ and $\beta$ are surjective we get
$$ \alpha \circ \epsilon = Id_{M}, \ \ \ \beta \circ \epsilon = Id_{M},$$
where $Id_{M}$ is the identity at $M$.\\
\noindent{Sometimes $M$ is denoted by $\Gamma_{\left(0\right)}$ and it is identified with the set $\epsilon \left(M\right)$ of identities of $\Gamma$. $\Gamma$ is also denoted by $\Gamma_{\left(1\right)}$. The elements of $M$ are called \textit{objects} and the elements of $\Gamma$ are called \textit{morphishms}. The map $\left(\alpha , \beta\right) : \Gamma \rightarrow M \times M$ is called the \textit{anchor map} and the space of sections of the anchor map is denoted by $\Gamma_{\left(\alpha, \beta\right)} \left(\Gamma\right)$. Finally, for each $g \in \Gamma$ the element $i \left( g \right)$ is denoted by $g^{-1}$.}\\
Roughly speaking, a groupoid may be depicted as a set of ``\textit{arrows}'' ($\Gamma$) joining points ($M$), in such a way that any two arrows may composed if the ending point of one coincides with the starting point of the other. Then, assuming natural conditions derived of the properties of a composition in a group, we get the definition of groupoid.

\begin{example}\label{5}
\rm
A group $G$ is a groupoid over a point and the operation law of the groupoid, $\cdot$, is the operation in $G$.
\end{example}




\noindent{Next, let us describe the crucial example of groupoid for the purpose of this paper.}
\begin{example}\label{8}
\rm
Let $A$ be a vector bundle on a manifold $M$. Denote by $A_{z}$, the fibre of $A$ over a $z \in M$. Then, the set $\Phi \left(A\right)$, consisting of all linear isomorphisms $L_{x,y}: A_{x} \rightarrow A_{y}$ for any $x,y \in M$, may be endowed with the structure of groupoid with structure maps,
\begin{itemize}
\item[(i)] $\alpha\left(L_{x,y}\right) = x$
\item[(ii)] $\beta\left(L_{x,y}\right) = y$
\item[(iii)] $L_{y,z} \cdot G_{x,y} = L_{y,z} \circ G_{x,y}, \ L_{y,z}: A_{y} \rightarrow A_{z}, \ G_{x,y}: A_{x} \rightarrow A_{y}$
\end{itemize}
We will call this groupoid as the \textit{frame groupoid on $A$}.\\
A relevant case is the \textit{1-jets groupoid on} $M$ and it arises when $A$ is the tangent bundle $TM$ of $M$. This groupoid is denoted by $\Pi^{1} \left(M,M\right)$. Notice that any isomorphism $L_{x,y}: T_{x}M \rightarrow T_{y}M$ may be written as a $1-$jet $j_{x,y}^{1} \psi$ of a local diffeomorphism $\psi$ from $M$ to $M$. Recall that the $1-$jet $j_{x,y}^{1} \psi$ may be identified with the tangent map $T_{x}\psi: T_{x}M \rightarrow T_{y}M$ (see \cite{SAUND} for details).
\end{example}

\begin{definition}\label{isotropygroup344}
\rm

A \textit{subgroupoid} of a groupoid $\Gamma \rightrightarrows M$ is a groupoid $\Gamma' \rightrightarrows M'$ such that $M' \subseteq M$, $\Gamma' \subseteq \Gamma$ and the the structure maps of $\Gamma'$ are the restrictions of the structure maps of $\Gamma$.
\end{definition}
Notice that the composition law of a subgroupoid is the same as that of the correspondent groupoid.
\begin{definition}\label{58}
\rm
Let $\Gamma \rightrightarrows M$ be a groupoid with $\alpha$ and $\beta$ the source map and target map, respectively. For each $x \in M$, the set
$$\Gamma^{x}_{x}= \beta^{-1}\left(x\right) \cap \alpha^{-1}\left(x\right),$$
is called the \textit{isotropy group of} $\Gamma$ at $x$. The set
$$\mathcal{O}\left(x\right) = \beta\left(\alpha^{-1}\left(x\right)\right) = \alpha\left(\beta^{-1}\left( x\right)\right),$$
is called the \textit{orbit} of $x$, or \textit{the orbit} of $\Gamma$ through $x$.
\end{definition}

\noindent{Observe that the isotropy groups inherit a \textit{bona fide} group structure.}
\begin{definition}
\rm
If $\mathcal{O}\left(x\right) = M$ for all $x \in M$ (or equivalently $\left(\alpha,\beta\right) : \Gamma  \rightarrow M \times M$ is a surjective map) the groupoid $\Gamma \rightrightarrows M$ is called \textit{transitive}. The sets,
$$ \alpha^{-1} \left(x \right) = \Gamma_{x}, \ \ \ \ \ \beta^{-1} \left(x \right) = \Gamma^{x},$$
are called $\alpha-$\textit{fibre at} $x$ and $\beta-$\textit{fibre at $x$}, respectively. We will denote
$$\Gamma_{x}^{y} = \Gamma_{x} \cap \Gamma^{y},$$
for all $x,y \in M$.
\end{definition}

\begin{definition}\label{9}
\rm
Let $\Gamma \rightrightarrows M$ be a groupoid. We may define the \textit{left translation by $g \in \Gamma$} as the map $L_{g} : \Gamma^{\alpha\left(g\right)} \rightarrow \Gamma^{\beta\left(g\right)}$, given by
$$ h \mapsto  g \cdot h .$$
We may define the \textit{right translation} by $g$, $R_{g} : \Gamma_{\beta\left(g\right)} \rightarrow \Gamma_{ \alpha \left(g\right)}$, analogously. 
\end{definition}
\noindent{Note that, the identity map on $\Gamma^{x}$ may be written as the following translation map,}
\begin{equation}\label{10} 
Id_{\Gamma^{x}} = L_{\epsilon \left(x\right)}.
\end{equation}
For any $ g \in \Gamma $, the left (resp. right) translation on $g$, $L_{g}$ (resp. $R_{g}$), is a bijective map with inverse $L_{g^{-1}}$ (resp. $R_{g^{-1}}$).\\\\

\begin{definition}
\rm
A \textit{Lie groupoid} is a groupoid $\Gamma \rightrightarrows M$ such that $\Gamma$ is a smooth manifold, $M$ is a smooth manifold and the structure maps are smooth. Furthermore, the source and the target map are submersions.\\
A \textit{Lie subgroupoid} of $\Gamma \rightrightarrows M$ is a Lie groupoid $\Gamma' \rightrightarrows M'$ such that it is a subgroupoid of $\Gamma$ satisfying that $\Gamma' $ and $M'$ are submanifolds of $\Gamma$ and $M$ respectively.
\end{definition}
As first example, any Lie group $G$ is a Lie groupoid (see example \ref{5}).

\begin{example}\label{15}
\rm
The frame groupoid $\Phi \left( A \right)$ on a vector bundle $A$ (see example \ref{8}) is a Lie groupoid . Let us consider two local coordinates, $\left(x^{i}\right)$ and $\left(y^{j}\right)$, on open neighbourhoods $U, V \subseteq M$, respectively, and two local basis of sections of $A_{U}$ and $A_{V}$, $\{\alpha_{p}\}$ and $\{ \beta_{q}\}$, respectively. The correspondent local coordinates $\left(x^{i} \circ \pi, \alpha^{p}\right)$ and $\left(y^{j} \circ \pi, \beta^{q}\right)$ on $A_{U}$ and $A_{V}$ are given by
\begin{itemize}
\item For any $a \in A_{U}$,
$$ a =  \alpha^{p}\left(a\right) \alpha_{p}\left(x^{i}\left(\pi \left(a\right)\right)\right).$$\\
\item For any $a \in A_{V}$,
$$ a =  \beta^{q}\left(a\right) \beta_{q}\left(y^{j}\left(\pi \left(a\right)\right)\right).$$
\end{itemize}
\noindent{Then, we can construct a local coordinate system on $\Phi \left(A\right)$
$$ \Phi \left(A_{U,V}\right) : \left(x^{i} , y^{j}_{i}, y^{j}_{i}\right),$$
where, $A_{U,V} = \alpha^{-1}\left(U\right) \cap \beta^{-1}\left(V\right)$ and for each $L_{x,y} \in \alpha^{-1}\left(x\right) \cap \beta^{-1}\left(y\right) \subseteq \alpha^{-1}\left(U\right) \cap \beta^{-1}\left(V\right)$, we have}
\begin{itemize}\label{16}
\item $x^{i} \left(L_{x,y}\right) = x^{i} \left(x\right)$.
\item $y^{j} \left(L_{x,y}\right) = y^{j} \left( y\right)$.
\item $y^{j}_{i}\left( L_{x,y}\right)  = A_{L_{x,y}}$, where $A_{L_{x,y}}$ is the associated matrix to the induced map of $L_{x,y}$ using the local coordinates $\left(x^{i} \circ \pi, \alpha^{p}\right)$ and $\left(y^{j} \circ \pi, \beta^{q}\right)$.
\end{itemize}
In the particular case of the $1-$jets groupoid on $M$, $\Pi^{1}\left(M,M\right)$, the local coordinates will be denoted as follows
\begin{equation}\label{17}
\Pi^{1}\left(U,V\right) : \left(x^{i} , y^{j}, y^{j}_{i}\right),
\end{equation}
where, for each $ j^{1}_{x,y} \psi \in \Pi^{1}\left(U,V\right)$
\begin{itemize}
\item $x^{i} \left(j^{1}_{x,y} \psi\right) = x^{i} \left(x\right)$.
\item $y^{j} \left(j^{1}_{x,y}\psi \right) = y^{j} \left( y\right)$.
\item $y^{j}_{i}\left( j^{1}_{x,y}\psi\right)  = \dfrac{\partial \left(y^{j}\circ \psi\right)}{\partial x^{i}_{| x} }$.
\end{itemize}

\end{example}

\begin{example}\label{13}
\rm
Let $\pi: P \rightarrow M$ be a principal bundle with structure group $G$. Denote by $\phi :  P  \times G \rightarrow P$ the action of $G$ on $P$.\\
Now, suppose that $\Gamma \rightrightarrows P$ is a Lie groupoid, with $\overline{\phi} :  \Gamma \times G \rightarrow \Gamma$ a free and proper action of $G$ on $\Gamma$ such that, for each $h \in G$, the pair $\left( \overline{\phi}_{h} , \phi_{h}\right)$ is an isomorphism of Lie groupoids, i.e.,
\begin{equation*}
\alpha_{\Gamma} \left( \overline{\phi}_{h} \left(g\right)\right) = \phi_{h} \left(\alpha_{\Gamma} \left(g \right)\right), \ \ \ \ \ \ \ \beta_{\Gamma} \left( \overline{\phi}_{h} \left(g\right)\right) = \phi_{h} \left(\beta_{\Gamma} \left(g \right)\right),
\end{equation*}
where $\alpha_{\Gamma}$ and $\beta_{\Gamma}$ are the source and the target map of $\Gamma \rightrightarrows P$, and preserves the composition, i.e.,
$$\overline{\phi}_{h} \left( g_{1} \cdot g_{2}\right) = \overline{\phi}_{h} \left(g_{1}\right) \cdot \overline{\phi}_{h} \left(g_{2}\right), \ \forall \left(g_{1} , g_{2} \right) \in \Gamma_{\left(2\right)}.$$
Then, we can construct a Lie groupoid $\Gamma / G \rightrightarrows M$ such that the source map, $\overline{\alpha}$, and the target map, $\overline{\beta}$, are given by
$$ \overline{\beta}\left([ g ]\right) = \pi \left( \beta_{\Gamma} \left(g\right)\right), \ \ \overline{\alpha}\left([g ] \right) = \pi \left( \alpha_{\Gamma}\left(g\right)\right),$$
for all $ g \in \Gamma$, and $[\cdot]$ denotes the equivalence class in the quotient space $\Gamma / G$. These kind of Lie groupoids are called \textit{quotient Lie groupoids by the action of a Lie group}.
\end{example}

Next, as an important example, we will introduce the \textit{second-order non-holonomic groupoid}.\\\\
Let $M$ be a manifold and $FM$ the frame bundle over $M$. So, we can consider the $1-$jets groupoid on $FM$, $\Pi^{1} \left(FM, FM\right) \rightrightarrows FM$.\\
Thus, we denote by $J^{1}\left(FM \right)$ the subset of $\Pi^{1} \left(FM,FM\right)$ given by the $1-$jets $j^{1}_{\overline{X},\overline{Y}} \Psi$ of local automorphism $\Psi$ of $FM$ such that
$$ \Psi \left( v \cdot g\right) = \Psi \left( v \right) \cdot g, \ \forall v \in Dom\left( \Psi \right), \ \forall g \in Gl \left( n , \mathbb{R} \right).$$
Let $\left(x^{i}\right)$ and $\left(y^{j}\right)$ be local coordinate systems over two open sets $U,V \subseteq M$, the induced coordinate systems over $FM$ are denoted by
$$FU: \left(x^{i} , x^{i}_{j}\right)$$
$$FV: \left(y^{j} , y^{j}_{i}\right).$$
Hence, we can construct induced coordinates over $\Pi^{1}\left(FM,FM\right)$
$$\Pi^{1}\left(FU,FV\right) = \left( \alpha , \beta \right)^{-1}\left(U , V\right) : \left(\left(x^{i},x^{i}_{j}\right),\left( y^{j} , y^{j}_{i}\right),y^{j}_{,i}, y^{j}_{,ik},y^{j}_{i,k} ,y^{j}_{i,kl}\right) ,$$
where for each $j^{1}_{\overline{X},\overline{Y}} \Psi \in \Pi^{1}\left(FU,FV\right)$, we have
\begin{itemize}
\item $x^{i} \left(j^{1}_{\overline{X},\overline{Y}} \Psi\right) = x^{i} \left(\overline{X}\right)$
\item $x^{i}_{j} \left(j^{1}_{\overline{X},\overline{Y}} \Psi\right) =x^{i}_{j} \left(\overline{X}\right)$
\item $y^{j} \left(j^{1}_{\overline{X},\overline{Y}} \Psi\right) = y^{j} \left( \Psi \left(\overline{X}\right)\right)$
\item $y^{j}_{i}\left( j^{1}_{\overline{X},\overline{Y}} \Psi\right)  = y^{j}_{i} \left( \Psi \left( \overline{X}\right)\right)$
\item $y^{j}_{,i} \left( j^{1}_{\overline{X},\overline{Y}} \Psi \right) =  \dfrac{\partial \left(y^{j}\circ \Psi\right)}{\partial x^{i}_{|\overline{X}} }$
\item $y^{j}_{,ik} \left( j^{1}_{\overline{X},\overline{Y}} \Psi \right) =  \dfrac{\partial \left(y^{j}\circ \Psi\right)}{\partial {x^{i}_{k}}_{|\overline{X}} }$
\item $y^{j}_{i,k} \left( j^{1}_{\overline{X},\overline{Y}} \Psi \right) =  \dfrac{\partial \left(y^{j}_{i}\circ \Psi\right)}{\partial {x^{k}}_{|\overline{X}} }$
\item $y^{j}_{i,kl} \left( j^{1}_{\overline{X},\overline{Y}} \Psi \right) =  \dfrac{\partial \left(y^{j}_{i}\circ \Psi\right)}{\partial {x^{k}_{l}}_{|\overline{X}} }$
\end{itemize}
Then, using these coordinates, $J^{1}\left(FM\right)$ can be described as follows:
{\footnotesize $$J^{1}\left(FU,FV\right) = J^{1}\left(FM\right) \cap \left(\alpha , \beta \right)^{-1}\left(U,V\right) : \left(\left(x^{i},x^{i}_{j}\right),\left( y^{j} , y^{j}_{i}\right),y^{j}_{,i}, 0,y^{j}_{i,k},y^{j}_{i,kl}\right) ,$$}
where
$$ y^{j}_{i,kl} = \left( \sum_{m}y^{j}_{m} \left( x^{-1} \right)^{m}_{k} \right) \delta_{l}^{i}.$$
Thus, $J^{1}\left(FM\right)$ is a submanifold of $\Pi^{1}\left(FM,FM\right)$ and its induced local coordinates will be denoted by
\begin{equation}\label{101}
J^{1}\left(FU,FV\right) : \left(\left(x^{i},x^{i}_{j}\right),\left( y^{j} , y^{j}_{i}\right),y^{j}_{,i},y^{j}_{i,k}\right) .
\end{equation}
Finally, restricting the structure maps we can ensure that $J^{1}\left(FM\right) \rightrightarrows FM$ is a Lie subgroupoid of the $1-$jets groupoid over $FM$.\\
We may now construct $j^{1}\left(FM\right)$ as the set of the $1-$jets of the form $j^{1}_{X,Y} F\psi$, where $\psi: M \rightarrow M$ is a local diffeomorphism. Let $\left( x^{i} \right)$ be a local coordinate system on $M$; then, restricting the induced local coordinates given in Eq. (\ref{101}) to $j^{1} \left( FM \right)$ we have that
$$ y^{j}_{i} = y^{j}_{,l}x^{l}_{i} \ \ \ \ \ ; \ \ \ \ \ y^{j}_{i,k} = y^{j}_{k,i}.$$
We deduce that $j^{1}\left(FM\right) \rightrightarrows FM$ is a reduced Lie subgroupoid of the $1-$jets groupoid over $FM$ and we denoted the coordinates on $j^{1} \left(FM \right)$ by
\begin{equation}\label{140}
j^{1}\left(FU,FV\right) : \left(\left(x^{i},x^{i}_{j}\right),\left( y^{j} , y^{j}_{i}\right),y^{j}_{i,k}\right) , \ \ \ y^{j}_{i,k} = y^{j}_{k,i}.
\end{equation}
Now, we will work with a quotient space of $J^{1}\left(FM\right)$ (resp. $j^{1}\left(FM\right)$) which will be our \textit{non-holonomic groupoid of second order} (resp. \textit{holonomic groupoid of second order}).\\\\
We consider the following right action of $Gl\left(n , \mathbb{R}\right)$ over $J^{1} \left(FM\right)$,

\begin{equation}\label{44}
\begin{array}{rccl}
\Phi : & J^{1} \left(FM\right) \times Gl \left( n , \mathbb{R}\right)  & \rightarrow & J^{1}\left(FM\right) \\
&\left(j^{1}_{\overline{X},\overline{Y}} \Psi , g\right)   &\mapsto &  j^{1}_{\overline{X} \cdot g , \overline{Y} \cdot g} \Psi.
\end{array}
\end{equation}

Thus, for each $g \in Gl \left( n , \mathbb{R}\right)$ the pair $\left(\Phi_{g} , R_{g}\right)$ (where $R$ is the natural right action of $Gl\left(n , \mathbb{R}\right)$ over $FM$) is a Lie groupoid automorphism. Therefore, we can consider the quotient Lie groupoid by this action $ \tilde{J}^{1} \left( FM \right) \rightrightarrows M$ which is called \textit{second-order non-holonomic groupoid over $M$}.\\
We will denote the structure maps of $\tilde{J}^{1} \left(FM\right)$ by $\overline{\alpha}$ and $\overline{\beta}$ (source and target maps respectively), $\overline{\epsilon}$ (identities map) and $\overline{i}$ (inversion map). The elements of $\tilde{J}^{1} \left( FM \right)$ are denoted by $j^{1}_{x,y} \Psi$ with $x,y \in M$ and $\overline{\alpha}\left( j^{1}_{x,y} \Psi\right) = x$ and $\overline{\beta}\left(j^{1}_{x,y} \Psi\right) = y$.\\
Then, the induced local coordinates are given by
\begin{equation}\label{45}
\tilde{J}^{1}\left(FU,FV\right) = \left( \overline{\alpha} , \overline{\beta} \right)^{-1}\left(U , V\right) : \left(\left(x^{i}\right),\left( y^{j} , y^{j}_{i}\right),y^{j}_{,i},y^{j}_{i,k}\right) .
\end{equation}
Considering ${e_{1}}_{x}$ as the $1-$jet through $x \in M$ which satisfies that $x^{i}_{j} \left( {e_{1}}_{x}\right) = \delta^{i}_{j}$ for all $i,j$, for each $j_{x,y}^{1}\Psi \in \tilde{J}^{1}\left(FM\right)$ we have
\begin{itemize}
\item $x^{i} \left(j^{1}_{x,y} \Psi\right) = x^{i} \left(x\right)$
\item $y^{j} \left(j^{1}_{x,y} \Psi\right) = y^{j} \left( y\right)$
\item $y^{j}_{i}\left( j^{1}_{x,y} \Psi\right)  = y^{j}_{i} \left( \Psi \left( {e_{1}}_{x}\right)\right)$.
\item $y^{j}_{,i} \left( j^{1}_{x,y} \Psi \right) =  \dfrac{\partial \left(y^{j}\circ \Psi\right)}{\partial x^{i}_{|x} }$
\item $y^{j}_{i,k} \left( j^{1}_{x,y} \Psi \right) =  \dfrac{\partial \left(y^{j}_{i}\circ \Psi\right)}{\partial {x^{k}}_{|{e_{1}}_{x}} }$
\end{itemize}

Observe that we can restrict the action $\Phi$ to an action of $Gl \left( n , \mathbb{R} \right)$ over $j^{1} \left(FM\right)$. So, by quotienting, we can build a reduced subgroupoid of $ \tilde{J}^{1} \left( FM \right) \rightrightarrows M$ which is denoted by $ \tilde{j}^{1} \left( FM \right) \rightrightarrows M$ and is called \textit{second-order holonomic groupoid over $M$}. Finally, by restriction, the local coordinates on $j^{1} \left( FM \right)$ are given by
\begin{equation}\label{secondholo345}
\tilde{j}^{1}\left(FU,FV\right) : \left(\left(x^{i}\right),\left( y^{j} , y^{j}_{i}\right),y^{j}_{i,k}\right) , \ \ \ y^{j}_{i,k} = y^{j}_{k,i}.
\end{equation}
Denote the structure maps of the holonomic groupoid over $M$ $ \tilde{j}^{1} \left( FM \right)$ by $\tilde{\alpha}$, $\tilde{\beta}$, $\tilde{\epsilon}$ and $\tilde{i}$.\\
Finally, let us define two projections $\overline{\Pi}_{1}^{2}$ and $\tilde{\Pi}_{1}^{2}$ from the non-holonomic groupoid $\tilde{J}^{1} \left(FM\right)$ of second order, to the $1-$jets groupoid $\Pi^{1}\left(M,M\right)$, as follows,
$$
\begin{array}{rccl}
\overline{\Pi}_{1}^{2} : & \tilde{J}^{1}\left(FM\right)  & \rightarrow & \Pi^{1}\left(M,M\right)\\
& j_{x,y}^{1}\Psi  &\mapsto &  \Psi\left(\overline{X}\right)[\overline{X}^{-1}]
\end{array}
$$
where $\overline{X} \in FM$ is a frame at $x$. It is easy to show that $\overline{\Pi}_{1}^{2}$ is well-defined and, locally,
$$\overline{\Pi}_{1}^{2} \left(\left(x^{i}\right),\left( y^{j} , y^{j}_{i}\right),y^{j}_{,i},y^{j}_{i,k}\right) = \left(x^{i}, y^{j} , y^{j}_{i}\right).$$
On the other hand we consider
$$
\begin{array}{rccl}
\tilde{\Pi}_{1}^{2} : & \tilde{J}^{1}\left(FM\right)  & \rightarrow & \Pi^{1}\left(M,M\right)\\
& j_{x,y}^{1}\Psi  &\mapsto &  j_{x,y}^{1}\psi
\end{array}
$$
where $\psi$ is the induced map of $\Psi$ over $M$. Then, locally
$$\tilde{\Pi}_{1}^{2} \left(\left(x^{i}\right),\left( y^{j} , y^{j}_{i}\right),y^{j}_{,i},y^{j}_{i,k}\right) = \left(x^{i}, y^{j} , y^{j}_{,i}\right).$$
Notice that $\overline{\Pi}_{1}^{2}$ and $\tilde{\Pi}_{1}^{2}$ are, indeed, Lie groupoid morphims over the identity map on $M$ (see \cite{KMG}).\\

Let us consider,
$$g = \left(\left(z^{i}\right),\left( y^{j} , g^{j}_{i}\right),g^{j}_{,i},g^{j}_{i,k}\right) , \ \ F= \left(\left(x^{i}\right),\left( z^{j} , F^{j}_{i}\right),F^{j}_{,i},F^{j}_{i,k}\right)\in \tilde{J}^{1} \left( F \mathcal{B} \right)$$
Then,
$$
 g \cdot F  = \left(\left(x^{i}\right),\left( y^{j} , \ g^{j}_{m}F^{m}_{i}\right),\ g^{j}_{,m}F^{m}_{,i},\ g^{j}_{r,m}F^{r}_{i}F^{m}_{,k} + g^{j}_{m}F^{m}_{i,k}\right)
$$
Hence, the left-translation by $g$ is given by,

\begin{equation}\label{lefttransl348}
\left(\left(x^{i}\right),\left( y^{j} , \ g^{j}_{m}y^{m}_{i}\right),\ g^{j}_{,m}y^{m}_{,i},\ g^{j}_{r,m}y^{r}_{i}y^{m}_{,k} + g^{j}_{m}y^{m}_{i,k}\right)
\end{equation}

Thus, the induced tangent map is characterized by the following equilities:
\begin{itemize}
    \item[i)] $$TL_{g} \left(\dfrac{\partial  }{\partial x^{i} }\right) = \dfrac{\partial  }{\partial x^{i} }$$
    \item[ii)] $$TL_{g} \left(\dfrac{\partial  }{\partial y^{j}_{i} }\right) = g^{m}_{j}\dfrac{\partial  }{\partial y^{m}_{i} } + g^{l}_{j,m}y^{m}_{,r} \dfrac{\partial  }{\partial y^{l}_{i,r} }$$
    \item[iii)] $$TL_{g} \left(\dfrac{\partial  }{\partial y^{j}_{,i} }\right) = g^{m}_{,j}\dfrac{\partial  }{\partial y^{m}_{,i} } + g^{l}_{n,j}y^{n}_{r} \dfrac{\partial  }{\partial y^{l}_{r,i} }$$
    
\item[iii)] $$TL_{g} \left(\dfrac{\partial  }{\partial y^{j}_{i,k} }\right) = g^{m}_{j}\dfrac{\partial  }{\partial y^{m}_{i,k} } $$
\end{itemize}
These equalities will be useful in what follows.\\\\

It could rise the situation in which we need to work with a (non necessarily Lie) subgroupoid of a Lie groupoid. To deal with this case, we have the so-called \textit{characteristic distribution} \cite{VMMDME,CHARDIST}.\\

\noindent{Let $ \Gamma \rightrightarrows M$ be a Lie groupoid and $\overline{\Gamma}$ be a subgroupoid of $\Gamma$ (not necessarily a Lie subgroupoid of $\Gamma$) over the same manifold $M$. We will denote by $\overline{\alpha}$, $\overline{\beta}$, $\overline{\epsilon}$ and $\overline{i}$ the restrictions of the structure maps $\alpha$, $\beta$, $\epsilon$ and $i$ of $\Gamma$ to $\overline{\Gamma}$ (see the diagram below)}\\

\begin{center}
 \begin{tikzcd}[column sep=huge,row sep=huge]
\overline{\Gamma}\arrow[r, hook, "j"] \arrow[rd, shift right=0.5ex] \arrow[rd, shift left=0.5ex]&\Gamma \arrow[d, shift right=0.5ex] \arrow[d, shift left=0.5ex] \\
& M 
 \end{tikzcd}
\end{center}

\noindent{where $j$ is the inclusion map. Now, we can construct a distribution $A \overline{\Gamma}^{T}$ over the manifold $\Gamma$ in the following way,
$$ g \in \Gamma \mapsto A \overline{\Gamma}^{T}_{g} \leq T_{g} \Gamma,$$
such that $A \overline{\Gamma}^{T}_{g}$ is the fibre of $ A \overline{\Gamma}^{T}$ at $g$ and it is generated by the (local) left-invariant vector fields $\Theta \in \frak X_{loc} \left( \Gamma \right)$ whose flow at the identities is totally contained in $\overline{\Gamma}$, i.e.,
\begin{itemize}
\item[(i)] $\Theta$ is tangent to the $\beta-$fibres, 
$$ \Theta \left( g \right) \in T_{g} \beta^{-1} \left( \beta \left( g \right) \right),$$
for all $g$ in the domain of $\Theta$.
\item[(ii)] $\Theta$ is invariant by left translations,
$$ \Theta \left( g \right) = T_{\epsilon \left( \alpha \left( g \right) \right) } L_{g} \left( \Theta \left( \epsilon \left( \alpha \left( g \right) \right) \right) \right),$$
for all $g $ in the domain of $\Theta$.
\item[(iii)] The (local) flow $\varphi^{\Theta}_{t}$ of $\Theta$ satisfies
$$\varphi^{\Theta}_{t} \left( \epsilon \left( x \right)\right) \in \overline{\Gamma}, $$
for all $x \in M$.
\end{itemize}
Notice that, for each $g \in \Gamma$, the zero vector $0_{g} \in T_{g} \Gamma$ is contained in the fibre of the distribution at $g$, namely $A \overline{\Gamma}^{T}_{g}$. On the other hand, it is easy to prove that a vector field $\Theta$ satisfies conditions (i) and (ii) if, and only if, its local flow $\varphi^{\Theta}_{t}$ is left-invariant or, equivalently,
$$ L_{g} \circ \varphi^{\Theta}_{t} = \varphi^{\Theta}_{t} \circ L_{g}, \ \forall g,t.$$
Then, taking into account that all the identities are in $\overline{\Gamma}$ (because it is a subgroupoid of $\Gamma$), condition (iii) is equivalent to the following,
\begin{itemize}
\item[(iii)'] The (local) flow $\varphi^{\Theta}_{t}$ of $\Theta$ at $\overline{g}$ is totally contained in $\overline{\Gamma}$, for all $\overline{g} \in \overline{\Gamma}$.
\end{itemize}
Thus, we are considering the \textit{left-invariant vector fields on $\Gamma$ whose integral curves are confined inside or outside $\overline{\Gamma}$}. It is also remarkable that, by construction, this distribution is differentiable, i.e., for each point $x$ and for any vector $v_{x}$ of the distribution at $x$ there exists a (local) vector field $\Theta$ tangent to the distribution such that,
$$ \Theta \left( x \right) = v_{x}.$$
\indent{The distribution $A \overline{\Gamma}^{T}$ is called the \textit{characteristic distribution of $\overline{\Gamma}$}.} For the sake of simplicity, we will denote the family of the vector fields which satisfy conditions (i), (ii) and (iii) by $\mathcal{C}$. The local vector fields of $\mathcal{C}$ will be called \textit{admissible vector fields for the couple $\left( \Gamma , \overline{\Gamma}\right)$}.\\
The structure of groupoid permits us to construct two more new objects associated to the distribution $A \overline{\Gamma}^{T}$. The first one is a smooth distribution over the base $M$ denoted by $A \overline{\Gamma}^{\sharp}$, called \textit{base-characteristic distribution}. The second one is a ``\textit{differentiable}'' correspondence $A\overline{\Gamma}$ which associates to any point $x$ of $M$ a vector subspace of $T_{\epsilon \left( x \right) } \Gamma$. Both constructions are characterized by the commutativity of the following diagram

\begin{large}
\begin{center}
 \begin{tikzcd}[column sep=huge,row sep=huge]
\Gamma\arrow[r, "A \overline{\Gamma}^{T}"] &\mathcal{P} \left( T \Gamma \right) \arrow[d, "T\alpha"] \\
 M \arrow[u,"\epsilon"] \arrow[r,"A \overline{\Gamma}^{\sharp}"] \arrow[ru,dashrightarrow, "A \overline{\Gamma}"]&\mathcal{P} \left( T M \right)
 \end{tikzcd}
\end{center}
\end{large}

%

\vspace{15pt}
\noindent{where $\mathcal{P} \left( E \right)$ defines the power set of $E$. Therefore, for each $x \in M$, the fibres satisfy that}
\begin{eqnarray*}
A \overline{\Gamma}_{x} &=&  A \overline{\Gamma}^{T}_{\epsilon \left( x \right)}\\
A \overline{\Gamma}^{\sharp}_{x}  &=& T_{\epsilon \left( x \right) } \alpha \left( A \overline{\Gamma}_{x} \right)
\end{eqnarray*}
It is remarkable that all the distributions introduced are not, necessarily, regular.\\
Notice that, taking into account that $A \overline{\Gamma}^{T}$ is locally generated by left-invariant vector field, we have that for each $g \in \Gamma$,
$$ A \overline{\Gamma}^{T}_{g} = T_{\epsilon \left( \alpha \left( g \right) \right)} L_{g} \left( A \overline{\Gamma}^{T}_{\epsilon \left( \alpha \left( g \right)\right)} \right),$$
i.e., the characteristic distribution is \textit{left-invariant}.\\

\vspace{0.15cm}

\indent{To summarize, associated to $\overline{\Gamma}$ we have three differentiable objects $A \overline{\Gamma}$, $A \overline{\Gamma}^{T}$ and $A \overline{\Gamma}^{\sharp}$.}
Now, we will study how these objects endow $\overline{\Gamma}$ with a sort of ``\textit{differentiable}'' structure.} In particular, by using the Stefan-Sussman's theorem \cite{HJS,PS} which deals with the integrability of singular distributions, we may prove the following result:

\begin{theorem}\label{10.24}
Let $\Gamma \rightrightarrows M$ be a Lie groupoid and $\overline{\Gamma}$ be a subgroupoid of $\Gamma$ (not necessarily a Lie groupoid) over $M$. Then, the characteristic distribution and the base characteristic distribution are integrable by the foliations $\overline{\mathcal{F}}$ and $\mathcal{F}$, respectively. Furthermore, $\overline{\Gamma}$ is a union of leaves of $\overline{\mathcal{F}}$.
\end{theorem}
$\overline{\mathcal{F}}$ and $\mathcal{F}$ are called \textit{characteristic foliation} and \textit{base-characteristic foliation}, respectively. Observe that,
\begin{itemize}
\item[(i)] For each $g \in \Gamma$,
$$\overline{\mathcal{F}} \left( g \right) \subseteq \Gamma^{\beta \left( g \right)}.$$
Indeed, if $g \in \overline{\Gamma}$, then
$$\overline{\mathcal{F}} \left( g \right) \subseteq \overline{\Gamma}^{\beta \left( g \right)}.$$
\item[(ii)] For each $g ,h \in \Gamma$ such that $\alpha \left( g \right) = \beta \left( h \right)$, we have
$$\overline{\mathcal{F}} \left( g \cdot h\right) = g \cdot \overline{\mathcal{F}} \left(  h\right).$$
\end{itemize}
It is remarkable that property $\left( i \right)$ means that each leaf of the foliation $\overline{\mathcal{F}}$ which integrates $A \overline{\Gamma}^{T}$ is contained in just one $\beta-$fibre, i.e., for each $g \in \Gamma$ the leaf $\overline{\mathcal{F}} \left( g \right)$ satisfies that 
$$\beta \left( h \right) = \beta \left( g \right),$$
for all $h \in \overline{\mathcal{F}} \left( g \right)$. Notice also that, one could expect that $\overline{\mathcal{F}} \left( g \right) = \overline{\Gamma}^{\beta \left( g \right)}$ but this is not true in general.\\
Observe that, any vector field $\Theta$ on $\Gamma$, may be projected onto a vector field $\Theta^{\sharp}$ on $M$ by the following identity:
$$ \Theta^{\sharp} \left( x \right) = T_{\epsilon\left(X\right)} \alpha \left( \Theta \left( \epsilon\left(X\right) \right)\right) , \ \forall x \in M$$
In fact, any admissible vector field $\Theta$ for the couple $\left( \Gamma , \overline{\Gamma}\right)$ projects into a vector field $\Theta^{\sharp}$ tangent to base-characteristic foliation $\mathcal{F}$.\\
As a complementary result, we may prove the following condition of maximality \cite{VMMDME,CHARDIST}.

\begin{proposition}\label{10.35}
Let $\overline{\mathcal{G}}$ be a foliation of $\Gamma$ such that $\overline{\Gamma}$ is a union of leaves of $\overline{\mathcal{G}}$ and 
$$\overline{\mathcal{G}} \left( g\right) \subseteq \Gamma^{ \beta \left( g \right)}, \ \forall g  \in \Gamma.$$
Then, the characteristic foliation $\overline{\mathcal{F}}$ is coarser that $\overline{\mathcal{G}}$, i.e.,
$$ \overline{\mathcal{G}} \left( g \right) \subseteq \overline{\mathcal{F}} \left( g \right)  , \ \forall g \in \Gamma.$$
\end{proposition}

Observe that, without assuming that $\overline{\Gamma}$ is a manifold, Theorem \ref{10.24} and Proposition \ref{10.35} prove that $\overline{\Gamma}$ may be divided into a \textit{maximal} union of leaves of a foliation of $\Gamma$, i.e., $\overline{\Gamma}$ may be divided into ``\textit{differentiable}'' parts in a maximal way. This gives us some kind of ``\textit{differentiable}'' structure over $\overline{\Gamma}$.\\

Let us now construct an algebraic structure of a groupoid over the leaves of $\mathcal{F}$. We will consider the minimal transitive groupoid $\overline{\Gamma} \left( \mathcal{F} \left( x \right) \right)$ generated by $\overline{\mathcal{F}} \left( \epsilon \left( x \right) \right) $. This groupoid is, in fact, generated by imposing that for all $\overline{g},\overline{h} \in \overline{\mathcal{F}} \left( \epsilon \left( x \right) \right) $ \cite{VMMDME},
$$\overline{g}, \overline{g}^{-1}, \overline{h}^{-1} \cdot \overline{g} \in \overline{\Gamma} \left( \mathcal{F} \left( x \right) \right).$$
Indeed, it is satisfies that,
\begin{equation}\label{10.17}
 \overline{\Gamma} \left( \mathcal{F} \left( x \right) \right) = \sqcup_{\overline{g} \in \overline{\mathcal{F}} \left( \epsilon \left( x \right) \right)} \overline{\mathcal{F}} \left( \epsilon \left( \alpha \left( \overline{g} \right) \right) \right),
\end{equation}
i.e., $\overline{\Gamma} \left( \mathcal{F} \left( x \right) \right)$ can be depicted as a disjoint union of fibres at the identities.\\
Observe that the $\beta-$fibre of this groupoid at a point $y \in \mathcal{F} \left( x \right)$ is given by $\overline{\mathcal{F}} \left( \epsilon \left( y \right) \right)$. Hence, the $\alpha-$fibre at $y$ is 
$$ \overline{\mathcal{F}}^{-1} \left( \epsilon \left( y\right) \right) = i \circ \overline{\mathcal{F}} \left( \epsilon \left( y \right) \right).$$
Furthermore, the Lie groups $\overline{\mathcal{F}} \left( \epsilon \left( y \right) \right) \cap \Gamma_{ y}$ are exactly the isotropy groups of $\overline{\Gamma}\left( \mathcal{F} \left( x \right) \right) $. All these results imply the following one (\cite{CHARDIST}):
\begin{theorem}\label{10.20}
For each $x \in M$ there exists a transitive Lie subgroupoid $\overline{\Gamma} \left( \mathcal{F} \left( x \right) \right)$ of $\Gamma$ with base $\mathcal{F} \left( x \right)$.
\end{theorem}

Thus, in fact, we have divided the manifold $M$ into leaves $\mathcal{F} \left( x \right)$ which have a maximal structure of transitive Lie subgroupoids of $\Gamma$.

As a particular consequence we have that: \textit{$\overline{\Gamma}$ is a transitive Lie subgroupoid of $\Gamma$ if, and only if, $M = \mathcal{F} \left( x \right)$ and $\overline{\Gamma} =\overline{\Gamma} \left( \mathcal{F} \left( x \right) \right)$ for some $x \in M$.}\\

\section{Cosserat Media}\label{CosseratSection234}

Here, we will give a very brief introduction of a model for materials with internal structure called \textit{Cosserat media} (\cite{FCosse,EPSTEIN1998127}).\\
Let us start with the so-called \textit{simple materials}. A \textit{body} $\mathcal{B}$ is modelized as a three-dimensional differentiable manifold and each point $X \in \mathcal{B}$ is called \textit{material particle} or \textit{material point}. The material points will be written using capital letters ($X,Y,Z, \dots$).\\
An embedding $\phi : \mathcal{B} \rightarrow \mathbb{R}^{3}$ is called a \textit{configuration of} $\mathcal{B}$. The $1-$jet $j_{X,\phi \left(X\right)}^{1} \phi$ of a configuration $\phi$ at $X \in \mathcal{B}$ is called an \textit{infinitesimal configuration at $X$}. The points $x \in \phi \left( \mathcal{B}\right) \subseteq \mathbb{R}^{3}$ are called \textit{spatial particles} or \textit{spatial points} and they will be denoted by lowercase letters ($x,y,z, \dots $).\\
We usually assume the existence of one of its configurations, say $\phi_{0}$, called \textit{reference configuration}, which is used to identify the body with an open subset of $\mathbb{R}^{3}$. Given any arbitrary configuration $\phi$, the change of configurations $\kappa = \phi \circ \phi_{0}^{-1}$ is called a \textit{deformation}, and its $1-$jet $j_{\phi_{0}\left(X\right) , \phi \left(X\right)}^{1} \kappa$ is called an \textit{infinitesimal deformation at $\phi_{0}\left(X\right)$}. Coordinates generated by the reference configuration will be denoted by capital letters $X^{I}$, while any other coordinates will be denoted by lowercase letters $x^{i}$.\\
For elastic simple bodies, the material is completely characterized by one function $W$ which depends, at each material particle $X \in \mathcal{B}$, on the gradient of the deformation evaluated at the point. Thus, $W$ is defined (see \cite{MELZA}) as a differentiable map
\begin{equation}\label{matresp34}
W : \mathcal{B} \times Gl\left(3 , \mathbb{R}\right) \rightarrow V,
\end{equation}
where $V$ is a real vector space. In general, $V$ will be the space of \textit{stress tensors} \cite{STRESS}. In fact, the contact forces at a particle $X$, in a given configuration $\phi$, are characterized by a symmetric second-order tensor $T_{X,\phi}$ on $\mathbb{R}^{3}$, which is called the stress tensor. Then, the mechanical response is given by the following equation:
$$W \left( X , F \right) = T_{X,\phi},$$
where $F$ is the $1-$jet at $\phi_{0} \left( X \right)$ of $\phi \circ \phi_{0}^{-1}$. Another equivalent way of considering $W$ is as a differentiable map
$$ W : \Pi^{1} \left( \mathcal{B}, \mathcal{B}\right) \rightarrow V,$$
by taking the associated matrix of the $1-$jets $j_{\phi_{0}\left( X\right),\phi_{0}\left( Y\right)}^{1}\left(\phi_{0}\circ \phi \circ \phi_{0}^{-1} \right)$ for each $j_{X,Y}^{1}\phi  \in  \Pi^{1} \left( \mathcal{B}, \mathcal{B}\right)$.\\

The picture describing the internal structure if a Cosserat medium is a little bit more complicated. In particular, the continuum will be endowed with the extra kinematic degrees of freedom provided by three independent deformable vectors attached at each material particle. So, a \textit{Cosserat medium} will be modelized by the linear frame bundle $F \mathcal{B}$ of a body $\mathcal{B}$. $\mathcal{B}$ is usually called the \textit{macromedium} or \textit{underlying body}. With some abuse of notation, we shall call $\mathcal{B}$ the \textit{Cosserat continuum}. Elements of $F\mathcal{B}$ will be denoted by overlined capital letters ($\overline{X}, \overline{Y}, \overline{Z}, \dots$) and the elements of $F\mathbb{R}^{3}$ will be denoted by overlined lowercase letters ($\overline{x}, \overline{y}, \overline{z}, \dots$).\\
A \textit{configuration} of a Cosserat medium $F\mathcal{B}$ is an embedding $\Psi : F \mathcal{B} \rightarrow F \mathbb{R}^{3}$ of principal bundles such that the induced Lie group morphism $\tilde{\psi} : Gl \left( 3 , \mathbb{R}\right) \rightarrow Gl\left(3, \mathbb{R}\right)$ is the identity map. Hence $\Psi$ satisfies
$$ \Psi \left( \overline{X} \cdot g \right) = \Psi \left(\overline{X}\right) \cdot g, \ \forall \overline{X} \in F \mathcal{B}, \ \forall g \in Gl \left( 3 , \mathbb{R} \right).$$
Notice that, $\Psi$ induces an embedding $\psi : \mathcal{B} \rightarrow \mathbb{R}^{3}$ verifying
$$ \pi_{\mathbb{R}^{3}} \circ \Psi = \psi \circ \pi_{\mathcal{B}}.$$
In particular, $\psi$ is an embedding of the macromedium $\mathcal{B}$ in $\mathbb{R}^{3}$. Furthermore, it satisfies that the subbundle $\Psi \left( F \mathcal{B}\right)$ of $F \mathbb{R}^{3}$ is just the frame bundle of $\psi \left( \mathcal{B}\right)$, i.e.,
$$\Psi \left( F \mathcal{B}\right) = F \psi \left( \mathcal{B}\right).$$
Since we are dealing with equivariants embedding, we can consider equivalence classes of the $1-$jets $j^{1}_{\overline{X} , \Psi \left( \overline{X}\right)} \Psi$ according to the action (\ref{44}). So, the equivalence class of an $1-$jet $j^{1}_{\overline{X} , \Psi \left( \overline{X}\right)} \Psi$, which is denoted by $j^{1}_{X,\psi \left(X\right)} \Psi$ like in the non-holonomic groupoid of second order, is called an \textit{infinitesimal configuration at $X$}. We usually identify the Cosserat medium with a configuration, say $\Psi_{0} : F \mathcal{B} \rightarrow F \mathbb{R}^{3}$, and we denote by $\psi_{0}$ the induced map of $\Psi_{0}$. $\Psi_{0}$ is called \textit{reference configuration}. Given any configuration $\Psi$, the change of configuration $\tilde{\kappa} = \Psi \circ \Psi_{0}^{-1}$ is called a \textit{deformation}, and its class of $1-$jets $j^{1}_{\psi_{0} \left(X\right) , \psi \left(X\right)} \tilde{\kappa}$ is called an \textit{infinitesimal deformation} at $\psi_{0}\left(X\right)$. Notice that the induced map of $\tilde{\kappa}$, is given by $\kappa = \psi \circ \psi_{0}^{-1}$.\\
From now on we make the following identification: $F \mathcal{B} \cong F \psi_{0} \left( \mathcal{B} \right)$.\\
Our assumption is that the material is completely characterized by one differentiable function $W : \tilde{J}^{1} \left( F \mathcal{B} \right) \rightarrow V$ over a vector space $V$. This map measures, for instance, the stored energy per unit mass and, again, we will call this function \textit{response functional} or \textit{mechanical response}. Notice that, by trivializing with the reference configuration, this map may be written as a differentiable map $W \left( X, \tilde{F}\right) $ depending on the particles $X$ on the macromedium and the Jacobian matrix $\tilde{F}$ of the (local) isomorphisms from $F \mathbb{R}^{3}$ to $F\mathbb{R}^{3}$.\\

Now, suppose that an infinitesimal neighbourhood of the material around the point $Y$ can be turned into a neighbourhood of $X$ such that the transformation cannot be detected by any mechanical experiment. If this condition is satisfied with every material particle $X$ of $ \mathcal{B}$, the body is said \textit{uniform}. We may express this physical property in a geometric way as follows.

\begin{definition}
\rm
A Cosserat continuum $\mathcal{B}$ is said to be \textit{uniform} if for each two points $X,Y \in \mathcal{B}$ there exists a local principal bundle isomorphism over the identity map on $Gl\left(3 , \mathbb{R}\right)$, $\Psi$, from $FU \subseteq F\mathcal{B}$ with $X\in U$ to $FV \subseteq F\mathcal{B}$ with $Y\in V$, where $U$ and $V$ are open neighbourhood of $M$, such that $\psi \left(X\right) =Y$ and
\begin{equation}\label{66}
W \left( j^{1}_{Y, \phi \left(Y\right)} \Phi \circ j^{1}_{X,Y} \Psi \right) = W \left( j^{1}_{Y, \phi \left(Y\right)} \Phi\right),
\end{equation}
for all infinitesimal configuration $j^{1}_{Y , \phi \left(Y\right)} \Phi$.
\end{definition}
$1-$jets $j^{1}_{X,Y} \Psi$ satisfying Eq. (\ref{66}) will be called \textit{material isomorphisms} from $X$ to $Y$ and will be relevant for the sequel. Two material points $X,Y$ will be called \textit{materially isomorphic} if there exists a material isomorphism from $X$ to $Y$. Notice that the relation of ``\textit{materially isomorphic}'' is an equivalence relation.\\
By composing the reference configuration with the material isomorphism, we may prove the following result:
\begin{proposition}\label{intuition}
Let $F\mathcal{B}$ be a Cosserat material. Two body particles $X$ and $Y$ are materially isomorphic if, and only if, there exist two (local) configurations $\Psi_{1}$ and $\Psi_{2}$ such that
$$ W_{1} \left( X , \tilde{F} \right)= W_{2} \left( Y , \tilde{F} \right), \ \forall \tilde{F},$$
where $W_{i}$ is the mechanical response associated to $\Psi_{i}$ for $i=1,2$, and 
\end{proposition}
This result provides us an intuition behind the notion of material isomorphism. In fact, two particles will be made of the same material if the mechanical response is the same under the action of two (possibly different) reference configurations.\\
For any two particles $X$ and $Y$, we will denote by $\overline{G} \left(X,Y\right)$ the collection of all $1-$jets $j_{X,Y}^{1}\Psi$ which satisfy Eq. (\ref{66}). So, the set $\overline{\Omega} \left( \mathcal{B}\right) = \cup_{X,Y \in \mathcal{B}} \overline{G}\left(X,Y\right)$ can be considered as a groupoid over $\mathcal{B}$ which is, indeed, a subgroupoid of the second-order non-holonomic groupoid $\tilde{J}^{1}\left(F \mathcal{B}\right)$. So, as an abuse of notation, we will denote the structure maps of $\overline{\Omega} \left( \mathcal{B} \right)$ using the same symbols used for the structure maps of $\tilde{J}^{1}\left(F \mathcal{B}\right)$. We will also denote $\overline{\alpha}^{-1} \left(X\right)$ (resp. $\overline{\beta}^{-1} \left(X\right)$) by $\overline{\Omega}_{X}\left( \mathcal{B}\right)$ (resp. $\overline{\Omega}^{X}\left( \mathcal{B}\right)$). $\overline{\Omega} \left(\mathcal{B}\right)$ is said to be the \textit{second-order non-holonomic material groupoid of $\mathcal{B}$}.\\

\begin{definition}\label{symm345}
\rm
Given a material point $X \in \mathcal{B}$ a \textit{material symmetry} at $X$ is a class of $1-$jets $j_{X,X}^{1}\Psi$, where $\Psi$ is a local automorphism at $X$ over the identity map on $Gl\left( 3 , \mathbb{R}\right)$, which satisfies Eq. (\ref{66}). 
\end{definition}
We denote by $\overline{G}\left(X\right)$ the set of all material symmetries which is, indeed, the isotropy group of $\overline{\Omega} \left(\mathcal{B}\right)$ at $X$ (see Definition \ref{isotropygroup344}). So, the following result is obvious.
\begin{proposition}
Let $\mathcal{B}$ be a Cosserat continuum. $\mathcal{B}$ is uniform if, and only if, $\overline{\Omega} \left( \mathcal{B}\right)$ is a transitive subgroupoid of $\tilde{J}^{1} \left( F \mathcal{B}\right)$.
\end{proposition}

Observe that the $\overline{\Omega} \left( \mathcal{B}\right)$ does not necessarily have a stucture of Lie groupoid. Indeed, notice that the definition of uniformity is a \textit{pointwise} property. In fact, consider a uniform Cosserat body $F\mathcal{B}$ and a fixed particle $X_{0}$, for any other particle $Y$ we may choose a material isomorphism from $Y$ to $X_{0}$, say $P \left( Y \right)$. So, we can construct a map $P : \mathcal{B} \rightarrow  \tilde{J}^{1} \left( F \mathcal{B}\right)$ consisting of material isomorphisms. Nevertheless, $P$ does not have to be differentiable. In other words, even when the Cosserat manifold is uniform, the choice of the material isomorphisms is not, necessarily, \textit{smooth}.
\begin{definition}\label{1.7.smoothuni}
\rm
A body $\mathcal{B}$ is said to be \textit{smoothly uniform} if for each material point $X \in \mathcal{B}$ there is a neighbourhood $\mathcal{U} $ around $X$ and a smooth map $P : \mathcal{U} \rightarrow \tilde{J}^{1} \left( F \mathcal{B}\right)$ such that for all $Y \in \mathcal{U}$ it satisfies that $P \left(Y\right)$ is a material isomorphism from $Y$ to $X$. The map $P$ is called a \textit{left (local) smooth field of material isomorphisms at} $X$. A \textit{right (local) smooth field of material isomorphisms at} $X$ will be a smooth map $P : \mathcal{U} \rightarrow \tilde{J}^{1} \left( F \mathcal{B}\right) $ such that for all $Y \in \mathcal{U}$ it satisfies that $P \left(Y\right)$ is a material isomorphism from $X$ to $Y$.
\end{definition}

Assume that $P$ is a right (local) smooth field of material isomorphisms at a material point $T$. Then, $P$ generates a smooth section $\mathcal{P}$ of the anchor map $\left( \overline{\alpha}, \overline{\beta}\right)$ of $\overline{\Omega}\left( \mathcal{B}\right)$ in the following way,
\begin{equation}\label{sections243}
 \mathcal{P}\left( X , Y \right) = P\left( Y\right)  \left[P\left( X \right)^{-1}\right]
 \end{equation}
The converse is also true. In other words, any smooth section $\mathcal{P}$ of the anchor map $\left( \overline{\alpha}, \overline{\beta}\right)$ of $\overline{\Omega}\left( \mathcal{B}\right)$, generates a right smooth field of material isomorphisms (and a left smooth field of material isomorphisms) satisfying Eq. (\ref{sections243}). Notice that, here, the word \textit{section} has a \textit{categorical meaning}; in fact, these sections should satisfy that
$$\mathcal{P}\left( Z , Y \right) \cdot \mathcal{P}\left( X , Z \right) = \mathcal{P}\left( X , Y \right), \ \forall X,Y,Z \in \mathcal{B}.$$

The smooth sections $\mathcal{P}$ of the anchor map $\left( \overline{\alpha}, \overline{\beta}\right)$ will be called \textit{smooth field of material isomorphisms}.\\
On the other hand, we may define a map
$$ \overline{W} \left( j^{1}_{T, \phi \left(T\right)} \Phi \right) = W \left( j^{1}_{T, \phi \left(T\right)} \Phi \right),$$
on the space of the $1-$jets of local diffeomorphisms $\Phi$ at a fixed material point $T$. We have that, for any $1-$jet $j^{1}_{Y, \psi \left(Y\right)} \Psi $
\begin{equation}\label{1.8}
\hspace{-0.3cm} W \left( j^{1}_{Y, \psi \left(Y\right)} \Psi  \right) = W\left( j^{1}_{Y, \psi \left(Y\right)} \Psi \circ \cdot  P \left( Y \right)\right) = \overline{W}\left( j^{1}_{Y, \psi \left(Y\right)} \Psi  \cdot  P \left( Y \right)\right)
\end{equation}
The meaning of Eq. (\ref{1.8}) is that the dependence of the mechanical response (near to a material particle) of the body coordinates is given by a multiplication of $j^{1}_{Y, \psi \left(Y\right)} \Psi$ to the right by a right smooth field of material isomorphisms.

So, we may prove the following result,
\begin{proposition}
Let $\mathcal{B}$ be a Cosserat continuum. $\mathcal{B}$ is smoothly uniform if, and only if, $\overline{\Omega} \left( \mathcal{B}\right)$ is a transitive Lie subgroupoid of $\tilde{J}^{1} \left( F \mathcal{B}\right)$.
\end{proposition}
In \cite{COSVME}, authors assume that $\overline{\Omega} \left( \mathcal{B}\right)$ is in fact a Lie subgroupoid to characterized properties like uniformity and homogeneity. Here, we will not assume this fact and, to deal with this problem, we refer to the \textit{characteristic distributions} \cite{CHARDIST,VMMDME}.\\

Let us present a particular case of this model, the so-called \textit{second-grade elastic materials} \cite{SecondMat}. In this case, it is assumed that the material response of the body $\mathcal{B}$ (see Eq. (\ref{matresp34})) not only depend on the first derivative of the configuration, but on both the first and the second gradient of the deformation. In other words, the mechanical response of a second-grade material is given by a differentiable map

$$W : \tilde{j}^{1} \left( F \mathcal{B} \right) \rightarrow V$$
over a vector space $V$, where $\tilde{j}^{1} \left( F \mathcal{B} \right)$ is the second-order holonomic groupoid over $\mathcal{B}$ (see Eq. (\ref{secondholo345})).\\
\begin{remark}
\rm

Let $F \mathcal{B}$ be a Cosserat medium, whose mechanical response is given by $W : \tilde{J}^{1} \left( F \mathcal{B} \right) \rightarrow V$. Then, taking into account that, the second-order holonomic groupoid $\tilde{j}^{1} \left( F \mathcal{B} \right)$ is a Lie subgroupoid of the second-order non-holonomic groupoid $\tilde{J}^{1} \left( F \mathcal{B} \right)$, we may restrict $W$ into a differentiable map $W : \tilde{j}^{1} \left( F \mathcal{B} \right) \rightarrow V$.\\
Thus, one could think that any Cosserat material may be studied as a second-grade material. However, due to the map $W$ encoded all the internal properties of the material, we may lose information with the restriction. So, the following question arise: How to ensure that the restriction does not induce a ``loss of information''. To answer this question, we refere to \cite{MAREMDL}.
\end{remark}

Let $\mathcal{B}$ be a Cosserat medium, with $W$ as mechanical response. According with the article \cite{MAREMDL}, we will say that $\mathcal{B}$ is a second-grade material, if all the material isomorphisms are natural prolongations to the frame bundle of the induced diffeomorphisms on the basis, i.e., all the material isomorphisms are $1-$jets $j_{X,Y}^{1} F \psi$, where $\psi$ is a local diffeomorphism on the body $\mathcal{B}$.\\

Let us consider the set $\Omega \left( \mathcal{B} \right)$ of $1-$jets $j^{1}_{X,Y} F \psi$ satisfying Eq. (\ref{66}), for a local diffeomorphism $\psi: \mathcal{B} \rightarrow \mathcal{B}$. In other words, $\Omega \left( \mathcal{B} \right)$  is the set of all material isomorphisms which are natural prolongations of local diffeomorphisms. It satisfies that $\Omega \left( \mathcal{B} \right)$ has the structure of subgroupoid of the second-order holonomic groupoid over $\mathcal{B}$,  $\tilde{j}^{1} \left( F \mathcal{B} \right)$. It is also true that $\Omega \left( \mathcal{B} \right)$ is a subgroupoid of the second-order non-holonomic material groupoid, $\overline{\Omega} \left(\mathcal{B}\right)$. In fact,

\begin{equation}
    \Omega \left( \mathcal{B} \right) \ = \  \overline{\Omega} \left(\mathcal{B}\right) \cap  \tilde{j}^{1} \left( F \mathcal{B} \right)
\end{equation}

$\Omega \left(\mathcal{B}\right)$ is called \textit{second-order holonomic material groupoid}. Therefore, for any Cosserat material, there always are two canonically defined groupoids, $\overline{\Omega} \left(\mathcal{B}\right)$ and $\Omega \left(\mathcal{B}\right)$, which are useful to study the constitutive properties of the material. In fact, we may use them to differenciate between second-grade material and \textit{bona-fide} Cosserat materials.\\

\begin{proposition}
Let $\mathcal{B}$ be a Cosserat medium, with $W$ as mechanical response. Then, $\mathcal{B}$ is a second-grade material if, and only if,
$$  \overline{\Omega} \left(\mathcal{B}\right)  \ = \ \Omega \left(\mathcal{B}\right)$$
\end{proposition}

\section{Cosserat characteristic distributions}\label{CosseratDistrib233}

Let $F \mathcal{B}$ be a Cosserat medium, whose mechanical response is given by $W : \tilde{J}^{1} \left( F \mathcal{B} \right) \rightarrow V$. Then, we have constructed two canonically defined groupoids, $\overline{\Omega} \left(\mathcal{B}\right)$ and $\Omega \left(\mathcal{B}\right)$, which are subgroupoids of the Lie groupoid $\tilde{J}^{1} \left( F \mathcal{B} \right)$. However, these two groupoids do not have to be Lie subgroupoids of $\tilde{J}^{1} \left( F \mathcal{B} \right)$ and, therefore, we are facing a framework in which the characteristic distributions may be constructed.\\

The characteristic distribution $A \overline{\Omega} \left( \mathcal{B} \right)^{T}$ of the non-holonomic material groupoid of second order will be called the \textit{non-holonomic material distribution of second-order}. On the other hand, the base-characteristic distribution $A \overline{\Omega} \left( \mathcal{B} \right)^{\sharp}$ will be called the \textit{non-holonomic body-material distribution of second order}.\\
Let $\Theta$ be an admissible vector field for the couple $\left(\tilde{J}^{1} \left( F \mathcal{B} \right), \overline{\Omega} \left( \mathcal{B} \right)\right)$, i.e., its (local) flow $\Psi^{\Theta}_{t}\left( \overline{\epsilon} \left( X \right) \right)$ at the identity $\overline{\epsilon}\left( X \right) = j_{X,X}^{1}Id_{F \mathcal{B}}$, where $Id_{F \mathcal{B}}$ is the identity over $F \mathcal{B}$, satisfies that
$$\Psi^{\Theta}_{t}\left( \overline{\epsilon} \left( X \right) \right) \subseteq \overline{\Omega} \left( \mathcal{B} \right)$$
for all $X \in \mathcal{B}$ and $t$ in the domain of the flow at $\overline{\epsilon} \left(X \right)$. Therefore, for any $g \in \tilde{J}^{1} \left( F \mathcal{B} \right)s$, we have
\begin{eqnarray*}
TW \left( \Theta \left( g \right)\right) &=&  \dfrac{\partial}{\partial t_{\vert 0}}\left(W \left(\Psi^{\Theta}_{t}\left( g \right) \right)\right) \\
&=&  \dfrac{\partial}{\partial t_{\vert 0}}\left(W \left(g \cdot\Psi^{\Theta}_{t}\left( \overline{\epsilon} \left( \overline{\alpha} \left( g \right) \right) \right) \right)\right)  \\
&=&  \dfrac{\partial}{\partial t_{\vert 0}}\left(W \left(g \right)\right)  = 0.
\end{eqnarray*}
Hence, we have that
\begin{equation}\label{14.22}
TW \left( \Theta \right) = 0
\end{equation}
The converse is proved in a similar way. Therefore, although the construction of the characteristic distribution is quite abstract, in the case of a Cosserat material this distribution may be completely described and calculated following Eq. (\ref{14.22}). In fact, let us consider a (local) left-invariant vector field $\Theta$ on $\tilde{J}^{1} \left( F \mathcal{B} \right)$. Therefore, we have 
$$ \Theta \ = \ \Theta^{i} \dfrac{\partial  }{\partial x^{i} } + \Theta^{j}_{i} \dfrac{\partial  }{\partial y^{j}_{i} } + \Theta^{j}_{,i} \dfrac{\partial  }{\partial y^{j}_{,i} } + \Theta^{j}_{i,k} \dfrac{\partial  }{\partial y^{j}_{i,k} }$$

Then, by using that $\Theta$ is left-invariant and Eq. (\ref{lefttransl348}), the local expression of $\Theta$ may be written as,
\begin{equation*}
 \Theta \ = \  \Theta^{i} \dfrac{\partial  }{\partial x^{i} } + y^{j}_{l}\Theta^{l}_{i}\dfrac{\partial  }{\partial y^{j}_{i} } +  y^{j}_{,l}\Theta^{l}_{,i}  \dfrac{\partial  }{\partial y^{j}_{,i} } + \left[y^{j}_{l,k}\Theta^{l}_{i}  + y^{j}_{i,l} \Theta^{l}_{,k} + y^{j}_{l} \Theta^{l}_{i,k} \right]\dfrac{\partial  }{\partial y^{j}_{i,k} }
\end{equation*}
where all the functions $\Theta^{i}$, $ \Theta^{i}_{j}$, $ \Theta^{i}_{,j}$, and $ \Theta^{i}_{j,k}$ depend on the material points of the body manifold $\mathcal{B}$ (which are given by $\alpha \left( g \right)$, when the vector field is evaluated on $g$). Thus, $\Theta$ is an admissible vector field for the couple $\left(\tilde{J}^{1} \left( F \mathcal{B} \right), \overline{\Omega} \left( \mathcal{B} \right)\right)$ if, and only if,
\begin{small}
\begin{equation}\label{materialeqcoss2345}
-\Theta^{i} \dfrac{\partial W }{\partial x^{i} } + \Theta^{l}_{i}\left[ y^{j}_{l}\dfrac{\partial  W}{\partial y^{j}_{i} } +  y^{j}_{l,k} \dfrac{\partial  W}{\partial y^{j}_{i,k} }\right] + \Theta^{l}_{,i} \left[ y^{j}_{,l}  \dfrac{\partial  W}{\partial y^{j}_{,i} } +  y^{j}_{m,l}\dfrac{\partial  W}{\partial y^{j}_{m,i} }\right] +   \Theta^{l}_{i,k} y^{j}_{l} \dfrac{\partial  W}{\partial y^{j}_{i,k} }= 0 .
\end{equation}

\end{small}
In other words, to construct the non-holonomic material distribution of second-order we have to find (local) functions on the body $\mathcal{B}$, $\Theta^{i}$, $ \Theta^{i}_{j}$, $ \Theta^{i}_{,j}$ and $ \Theta^{i}_{j,k}$ solving the linear equation (\ref{materialeqcoss2345}). In this way, the local functions $\Theta^{i}$ generate non-holonomic body-material distribution of second order $A \overline{\Omega}\left( \mathcal{B} \right)^{\sharp}$. This equation will be called  \textbf{non-holonomic material equation for Cosserat media}.\\

Let us now work with the holonomic material groupoid $\Omega \left( \mathcal{B}\right)$. First, the characteristic distribution $A \Omega \left( \mathcal{B} \right)^{T}$ of the holonomic material groupoid of second order will be called the \textit{holonomic material distribution of second-order}. On the other hand, the base-characteristic distribution $A \Omega \left( \mathcal{B} \right)^{\sharp}$ will be called the \textit{holonomic body-material distribution of second order}.\\

Next, let $\Theta$ be an admissible vector field for the couple $\left(\tilde{J}^{1} \left( F \mathcal{B} \right), \Omega \left( \mathcal{B} \right)\right)$. Then, its (local) flow $\Psi^{\Theta}_{t}\left( \overline{\epsilon} \left( X \right) \right)$ at the identity $\overline{\epsilon}\left( X \right) = j_{X,X}^{1}Id_{F \mathcal{B}}$ is totally contained in $\Omega \left( \mathcal{B} \right) \subseteq \tilde{j}^{1} \left( F \mathcal{B} \right)$. Therefore, this flow should be given by 
$$ \Psi^{\Theta}_{t}\left( \overline{\epsilon} \left( X \right) \right) = j_{{\psi^{\Theta}_{-t}} \left( X \right) , X}^{1} F\psi^{\Theta}_{t}$$
i.e., the flow of $\Theta$ is totally characterized by the flow of the projected vector field $\Theta^{\sharp} = T \overline{\alpha} \circ \Theta \circ \overline{\epsilon}$ ($\psi^{\Theta}_{t}$). This kind of vector fields are sometimes called \textit{complete lift of $\Theta^{\sharp}$}. So,
\begin{small}
$$ \Theta \ = \  -\Theta^{i} \dfrac{\partial  }{\partial x^{i} } + y^{j}_{l}\dfrac{\partial  \Theta^{l}}{\partial x^{i} }\dfrac{\partial  }{\partial y^{j}_{i} }+ y^{j}_{,l}\dfrac{\partial  \Theta^{l}}{\partial x^{i} }\dfrac{\partial  }{\partial y^{j}_{,i} } + \left[ y^{j}_{l,k} \dfrac{\partial  \Theta^{l}}{\partial x^{i} } + y^{j}_{i,l} \dfrac{\partial  \Theta^{l}}{\partial x^{k} } + y^{j}_{l} \dfrac{\partial^{2} \Theta^{l} }{\partial x^{i}\partial x^{k} } \right]\dfrac{\partial  }{\partial y^{j}_{i,k} },$$
\end{small}
where the functions $\Theta^{i}$ depends on the material points of the body manifold $\mathcal{B}$. Thus, $\Theta$ is an admissible vector field for the couple $\left(\tilde{J}^{1} \left( F \mathcal{B} \right), \Omega \left( \mathcal{B} \right)\right)$ if, and only if,

\begin{small}
\begin{equation}\label{holonomicmaterialeqcoss2245}
-\Theta^{i} \dfrac{\partial W }{\partial x^{i} } + \dfrac{\partial  \Theta^{l}}{\partial x^{i} }\left[ y^{j}_{l}\dfrac{\partial  W}{\partial y^{j}_{i} }  +  y^{j}_{,l}  \dfrac{\partial  W}{\partial y^{j}_{,i} } + y^{j}_{l,k} \dfrac{\partial  W}{\partial y^{j}_{i,k} } +  y^{j}_{m,l}\dfrac{\partial  W}{\partial y^{j}_{m,i} }\right] +    \dfrac{\partial^{2} \Theta^{l} }{\partial x^{i}\partial x^{k} } y^{j}_{l} \dfrac{\partial  W}{\partial y^{j}_{i,k} }= 0 .
\end{equation}

\end{small}

Therefore, to construct the holonomic material distribution of second order, we have to solve a second order partial differential equation (\ref{holonomicmaterialeqcoss2245}). In particular, the functions $\Theta^{i}$ solving the PDE (\ref{holonomicmaterialeqcoss2245}) generate the holonomic body-material distribution of second order. This equation will be called  \textbf{holonomic material equation for Cosserat media}.\\

The foliations associated with the (non-)holonomic material distribution of second order and the (non-)holonomic body-material distribution of second order will be called \textit{(non-)holonomic material foliation of second order} and \textit{(non-)holonomic body-material foliation of second order}, and denoted by $\left(\overline{\mathcal{NF}}\right) \ \overline{\mathcal{F}} $, $\left(\mathcal{NF}\right) \ \mathcal{F}$, respectively.\\
For each $X \in \mathcal{B}$, we will denote the Lie groupoids $\overline{\Omega} \left( \mathcal{B} \right)\left(\mathcal{NF}\left( X \right)\right)$ and $\Omega \left( \mathcal{B} \right)\left(\mathcal{F}\left( X \right)\right)$ by $\overline{\Omega} \left( \mathcal{NF} \left( X \right) \right)$ and $\Omega \left(\mathcal{F}\left( X \right)\right)$, respectively (see Theorem \ref{10.20}). Recall that $\Omega \left( \mathcal{B} \right)\left(\mathcal{NF}\left( X \right)\right)$ is a subgroupoid of $\overline{\Omega} \left( \mathcal{B} \right)\left(\mathcal{F}\left( X \right)\right)$.\\
Observe that in continuum mechanics a \textit{sub-body} of a body $\mathcal{B}$ is given by an open submanifold of $\mathcal{B}$.  Here, however, the foliation $\mathcal{F}$ gives us submanifolds of different dimensions (not only dimension 3). Thus, we will follow \cite{CHARDIST,GENHOM} for a more general definition:
\begin{definition}\label{materialsubm45623}
\rm
Let us consider a submanifold $\mathcal{P}$ of $\mathcal{B}$. Then, a \textit{Cosserat submanifold of $F\mathcal{B}$} is given by all the elements of $F \mathcal{B}$ at points of $\mathcal{P}$, which is denoted by $F \mathcal{P}$. In cases where it causes no confusion we will refer
to the Cosserat material as the submanifold $\mathcal{P}$.
\end{definition}
It is important to note that any Cosserat submanifold $\mathcal{P}$ inherits certain material structure from $\mathcal{B}$. In particular, the material response of a material submanifold $\mathcal{P}$ is measured by restricting $W$ to the $1-$jets of local isomorphisms $\Psi$ on $F\mathcal{B}$ from $F\mathcal{P}$ to $F\mathcal{P}$. However, it is easy to observe that a material submanifold of a Cosserat medium is not exactly a Cosserat material (dimension of $\mathcal{P}$ is not restricted to be three).\\

Related with \cite{MAREMDL}, we will say that a Cosserat submanifold $\mathcal{P}$ is a \textit{second-grade material submanifold}, if all the material isomorphisms from points on $\mathcal{P}$ to points on $\mathcal{P}$ are natural prolongations to the frame bundle of the induced diffeomorphisms on $\mathcal{B}$. In other words, all the material isomorphisms from $\mathcal{P}$ to $\mathcal{P}$ are $1-$jets $j_{X,Y}^{1}F \psi$, where $\psi$ is a local diffeomorphism on the body $\mathcal{B}$.\\

As a corollary of Theorem \ref{10.24} and Proposition \ref{10.35}, we have the following result.
\begin{theorem}\label{14.1}
The non-holonomic body-material foliation $\mathcal{NF}$ (resp. holonomic body-material foliation $\mathcal{F}$) divides the body $\mathcal{B}$ into maximal smoothly uniform Cosserat submanifolds (resp. second-grade material submanifolds).
\end{theorem}

It should be also observed that, in this case, ``maximal'' means that any other foliation $\mathcal{H}$ by smoothly uniform material submanifolds (resp. second-grade material submanifolds) is thinner than $\mathcal{NF}$ (resp. $\mathcal{F}$), i.e.,
$$ \mathcal{H} \left( X \right) \subseteq  \mathcal{NF} \left( X \right) \left(\text{resp. } \mathcal{F}\left( X \right) \right)  , \ \forall X \in  \mathcal{B}.$$
We should notice that this result provides us two different and intuitive divisions of the Cosserat material. First, a Cosserat material could be \textit{strictly non-uniform}. However, it may be maximally decomposed into ``\textit{(smoothly) uniform parts}'' and this decomposition is, in fact, a foliation $\mathcal{NF}$ of the macromedium.

\begin{theorem}\label{uniformityth}
Let $F\mathcal{B}$ be a Cosserat material whose mechanical response is denoted by $W$. Then, $\mathcal{B}$ is smoothly uniform if, and only if, the non-holonomic material equation for Cosserat media (\ref{materialeqcoss2345})
\begin{equation*}
 \Theta^{i} \dfrac{\partial W }{\partial x^{i} } + \Theta^{l}_{i}\left[ y^{j}_{l}\dfrac{\partial  W}{\partial y^{j}_{i} } +  y^{j}_{l,k} \dfrac{\partial  W}{\partial y^{j}_{i,k} }\right] + \Theta^{l}_{,i} \left[ y^{j}_{,l}  \dfrac{\partial  W}{\partial y^{j}_{,i} } +  y^{j}_{m,l}\dfrac{\partial  W}{\partial y^{j}_{m,i} }\right] +   \Theta^{l}_{i,k} y^{j}_{l} \dfrac{\partial  W}{\partial y^{j}_{i,k} }= 0 .
\end{equation*}
may be solved for any initial condition for the triple $\left(\Theta^{1} \left( X \right), \Theta^{2} \left( X \right), \Theta^{3} \left( X \right) \right)$, for all $X \in \mathcal{B}$.
\end{theorem}

Notice that, for any (local) admissible vector field 
$$ \Theta \ = \ \Theta^{i} \dfrac{\partial  }{\partial x^{i} } + \Theta^{j}_{i} \dfrac{\partial  }{\partial y^{j}_{i} } +  \Theta^{j}_{,i} \dfrac{\partial  }{\partial y^{j}_{,i} } +  \Theta^{j}_{i,k} \dfrac{\partial  }{\partial y^{j}_{i,k} }$$
it satisfies that, locally
$$ \Theta^{\sharp} \ = \ \Theta^{i} \dfrac{\partial  }{\partial x^{i}}$$
Therefore, the functions $\Theta^{i}$ solving Eq. (\ref{materialeqcoss2345}) correspond with the coordinates of the projection $\Theta^{\sharp}$ of the admissible vector fields $\Theta$ for the couple $\left(\tilde{J}^{1} \left( F \mathcal{B} \right), \overline{\Omega} \left( \mathcal{B} \right)\right)$. So, for each material particle $X \in \mathcal{B}$, we may consider the space
$$ \mathcal{C}_{X}^{\sharp} := \{ \Theta^{\sharp}\left( X \right) \ : \ \Theta \text{ is an admissible vector field}\}$$
Observe that the equation (\ref{materialeqcoss2345}) is linear with respect to the solutions $\left( \Theta^{i}, \Theta^{j}_{i}\right)$. Therefore, $ \mathcal{C}_{X}^{\sharp}$ is a vector subspace of $T_{X}\mathcal{B}$ and it is equal to the fibre of the non-holonomic material distribution of second-order at $X$, i.e.,
$$ \mathcal{C}_{X}^{\sharp} = A \overline{\Omega} \left( \mathcal{B}\right)^{\sharp}_{X}$$
Hence, we may reformulate Theorem \ref{uniformityth} in the following way:
\begin{theorem}
Let $F\mathcal{B}$ be a Cosserat material whose mechanical response is denoted by $W$. Then, $\mathcal{B}$ is smoothly uniform if, and only if, $A \overline{\Omega}\left( \mathcal{B}\right)^{\sharp}_{X}$ has dimension $3$ for all particle $X\in \mathcal{B}$.
\end{theorem}
Thus, among other conclusions, we only have to solve Eq. (\ref{materialeqcoss2345}) for initial conditions on a basis of $\mathbb{R}^{3}$.\\
In this way, as a summary, Eq. (\ref{materialeqcoss2345}) works to study the \textit{uniformity property of the material}. On the other hand, the holonomic material equation for Cosserat media (\ref{holonomicmaterialeqcoss2245}) will be useful to study both, the uniformity and the property of ``\textit{being a second-grade material}'', at the same time.\\
So, as a second division of the material $\mathcal{B}$ provided by Theorem \ref{14.1}, the Cosserat material may be maximally decomposed into ``\textit{(smoothly) uniform second-grade submanifolds}'' and this decomposition is, again, a foliation $\mathcal{F}$ of the macromedium. So, a \textit{smoothly uniform material is a second grade material if, and only if, the foliation $\mathcal{F}$ consists of only one leaf (equal to $\mathcal{B}$}).\\

\begin{theorem}\label{uniformityandsecondth}
Let $F\mathcal{B}$ be a Cosserat material whose mechanical response is denoted by $W$. Then, $\mathcal{B}$ is a smoothly uniform second-grade material if, and only if, the holonomic material equation for Cosserat media (\ref{holonomicmaterialeqcoss2245})
\begin{equation*}
-\Theta^{i} \dfrac{\partial W }{\partial x^{i} } + \dfrac{\partial  \Theta^{l}}{\partial x^{i} }\left[ y^{j}_{l}\dfrac{\partial  W}{\partial y^{j}_{i} }  +  y^{j}_{,l}  \dfrac{\partial  W}{\partial y^{j}_{,i} } + y^{j}_{l,k} \dfrac{\partial  W}{\partial y^{j}_{i,k} } +  y^{j}_{m,l}\dfrac{\partial  W}{\partial y^{j}_{m,i} }\right] +    \dfrac{\partial^{2} \Theta^{l} }{\partial x^{i}\partial x^{k} } y^{j}_{l} \dfrac{\partial  W}{\partial y^{j}_{i,k} }= 0 .
\end{equation*}
may be solved for any initial condition for the triple $\left(\Theta^{1} \left( X \right), \Theta^{2} \left( X \right), \Theta^{3} \left( X \right) \right)$, for all $X \in \mathcal{B}$.
\end{theorem}

Again, the functions $\Theta^{i}$ solving Eq. (\ref{holonomicmaterialeqcoss2245}) correspond with the coordinates of the projection $\Theta^{\sharp}$ of the admissible vector fields $\Theta$ for the couple $\left(\tilde{J}^{1} \left( F \mathcal{B} \right), \Omega \left( \mathcal{B} \right)\right)$ and the space of the evaluations of these vector fields onto a particle $X \in \mathcal{B}$ is the fibre $ A \Omega \left( \mathcal{B}\right)^{\sharp}_{X}$. Therefore, we may reformulate Theorem \ref{uniformityandsecondth} in the following way:
\begin{theorem}
Let $F\mathcal{B}$ be a Cosserat material whose mechanical response is denoted by $W$. Then, $\mathcal{B}$ is a smoothly uniform second-grade material if, and only if, $A \Omega \left( \mathcal{B}\right)^{\sharp}_{X}$ has dimension $3$ for all particle $X\in \mathcal{B}$.
\end{theorem}
In particular, a \textit{smoothly uniform Cosserat material is a second-grade material} if, and only if, the space of solution of Eq. (\ref{holonomicmaterialeqcoss2245}) has dimension $3$ at all the points. However, this result is not enough to characterize the \textit{second-grade character} of arbitrary (uniform or not) Cosserat medium.


\begin{theorem}\label{uniformityandsecondth23}
Let $F\mathcal{B}$ be a Cosserat material whose mechanical response is denoted by $W$. Then, all the material submanifolds of $\mathcal{B}$ given by its maximal division in smoothly uniform materials are second-grade material submanifolds, if and only if, at each particle $X$, all the (local) solutions $\Theta^{i}$ for Eq. (\ref{holonomicmaterialeqcoss2245}) generate all the solutions of Eq. (\ref{materialeqcoss2345}) by the following equalities $\Theta^{j}_{i} = \dfrac{\partial   \Theta^{j}}{\partial x^{i}} $, $\Theta^{j}_{,i} = \dfrac{\partial   \Theta^{j}}{\partial x^{i}} $ and $\Theta^{j}_{i,k} = \dfrac{\partial^{2}   \Theta^{j}}{\partial x^{i}\partial x^{k}} $.

\begin{proof}
Assume that, for each particle $X$, all the (local) solutions $\Theta^{i}$ for Eq. (\ref{holonomicmaterialeqcoss2245}) generate all the solutions of Eq. (\ref{materialeqcoss2345}). Then, equivalently,
$$ A \overline{\Omega}\left( \mathcal{B}\right)^{T}_{\overline{\epsilon}\left(X \right)} \ = \ A \Omega\left( \mathcal{B}\right)^{T}_{\overline{\epsilon}\left(X \right)},$$
for all $X \in \mathcal{B}$. Thus,

$$ \overline{\mathcal{NF}}\left( \overline{\epsilon}\left(X \right)\right) \ = \ \overline{\mathcal{F}}\left( \overline{\epsilon}\left(X \right)\right) , \ \forall X \in \mathcal{B}$$
Therefore, by construction we have
$$ \overline{\Omega}\left( \mathcal{NF}\left(X \right)\right) \ = \  \Omega\left( \mathcal{F}\left(X \right)\right),$$
for all $X \in \mathcal{B}$. The converse is is proved following a similar argument.
\end{proof}
\end{theorem}
In particular, assume that the non-holonomic material groupoid of second order $\overline{\Omega}\left(\mathcal{B}\right) $ is a Lie groupoid. Then, for all particle $X \in \mathcal{B}$,
$$ \overline{\beta}^{-1}\left( X \right) \ = \ \overline{\mathcal{NF}}\left( \overline{\epsilon} \left( X \right) \right)$$
Thus, in the conditions of Theorem \ref{uniformityandsecondth23}, one deduces,
$$ \overline{\beta}^{-1}\left( X \right) \ = \ \overline{\mathcal{F}}\left( \tilde{\epsilon} \left( X \right) \right)$$
Therefore, $\overline{\Omega}\left(\mathcal{B}\right) = \Omega\left(\mathcal{B}\right) $.

\begin{theorem}\label{uniformityandsecondth23456}
Let $F\mathcal{B}$ be a Cosserat material whose mechanical response is denoted by $W$ in such a way that the non-holonomic material groupoid of second order $\overline{\Omega}\left(\mathcal{B}\right) $ is a Lie groupoid. Then, $\mathcal{B}$ is a second-grade material, if and only if, at each particle $X$, all the (local) solutions $\Theta^{i}$ for Eq. (\ref{holonomicmaterialeqcoss2245}) generate all the solutions of Eq. (\ref{materialeqcoss2345}) by the following equalities $\Theta^{j}_{i} = \dfrac{\partial   \Theta^{j}}{\partial x^{i}} $, $\Theta^{j}_{,i} = \dfrac{\partial   \Theta^{j}}{\partial x^{i}} $ and $\Theta^{j}_{i,k} = \dfrac{\partial^{2}   \Theta^{j}}{\partial x^{i}\partial x^{k}} $.
\end{theorem}

Thus, the equations (\ref{materialeqcoss2345}) and (\ref{holonomicmaterialeqcoss2245}) are also useful to investigate when a Cosserat material (smoothly uniform or not) is, in fact, a second-grade material by comparing the solutions of both equations.

\section{Homogeneity}\label{Homogeneitysection24}
As we already know, a Cosserat medium is (smoothly) uniform if the function $W$ does depend on the point $X$ in a multiplicative way (Eq. (\ref{1.8})). In addition, a Cosserat continuum is said to be \textit{homogeneous} if we can choose a global section of the non-holonomic material groupoid of second order which is constant on the body, more precisely:

\begin{definition}\label{151}
\rm
A Cosserat medium $\mathcal{B}$ is said to be \textit{homogeneous} if it admits a global configuration $\Psi$ which induces a global section of $\left(\overline{\alpha} , \overline{\beta}\right)$ in $\overline{\Omega} \left( \mathcal{B}\right)$, $\mathcal{P}$, i.e., for each $X , Y \in \mathcal{B}$
\begin{equation}\label{homog234}
\mathcal{P}\left(X,Y\right) = j^{1}_{X,Y} \left(\Psi^{-1} \circ F\tau_{\psi\left(Y\right) - \psi \left(X\right)} \circ \Psi\right),
\end{equation} 
where $\tau_{\psi\left(Y\right) - \psi \left(X\right)}: \mathbb{R}^{3} \rightarrow \mathbb{R}^{3}$ denotes the translation on $\mathbb{R}^{3}$ by the vector $\psi\left(Y\right) - \psi \left(X\right)$, and $\psi$ is the induced map of $\Psi$ over the macromedium $\mathcal{B}$. $\mathcal{B}$ is said to be \textit{locally homogeneous} if there exists a covering of $\mathcal{B}$ by homogeneous open sets. From now on, in cases where it causes no confusion we will refer to local homogeneity as \textit{homogeneity}.
\end{definition}

\begin{proposition}
Let $F\mathcal{B}$ be a Cosserat medium. Then, $F\mathcal{B}$ is (locally) homogeneous if, and only if, there exist (local) reference configurations such that for the associated constitutive laws $W$ does not depend on the base points, i.e.
\begin{equation*}
W \left( j_{X,Y}^{1}\Phi_{1} \right) = W \left( j_{Z,T}^{1}\Phi_{2} \right),
\end{equation*}
whenever it satisfies the associated Jacobian matrix of $\Phi_{1}$ and $\Phi_{2}$ at $X$ and $Z$, respectively, are equal. 
\end{proposition}
Therefore, a material body is homogeneous if there exists a configuration such that the material response does not depend on the body points.

Notice that local homogeneity is obviously more restrictive than smooth uniformity. In fact, a homogeneous Cosserat body is a smoothly uniform body in which the (local) smooth fields of material isomorphisms (see Definition \ref{1.7.smoothuni}) may be chosen to be induced by configurations in the sense of Eq. (\ref{homog234}). Sections of $\left(\overline{\alpha} , \overline{\beta}\right)$ in $\overline{\Omega} \left( \mathcal{B}\right)$ given by Eq. (\ref{homog234}) will be called \textit{homogeneous sections}.\\
However, in a purely intuitive picture, homogeneity can be interpreted as the absence of defects. Thus, it would make sense to have a proper definition of homogeneity for non-uniform Cosserat media. In the literature we can already find some partial answer of this question (\cite{FGMA,FGM2} for FGM's, \cite{EPST,MGEOEPS} for laminated and bundle materials and \cite{GENHOM} for simple materials).
\begin{definition}\label{homogensubm}
Let $F\mathcal{B}$ be a Cosserat material and $\mathcal{P}$ be a submanifold of $\mathcal{B}$. $\mathcal{P}$ is said to be \textit{homogeneous} if, and only if, for all point $X \in \mathcal{P}$ there exists a local configuration $\Psi$ of $F\mathcal{B}$ on an open subset $ FU \subseteq F\mathcal{B}$, with $\mathcal{P} \subseteq U$, which satisfies that

$$ j_{Y,Z}^{1} \left( \Psi^{-1} \circ \tau_{ \psi\left( Z \right) - \psi \left( Y\right)} \circ \Psi \right),$$
is a material isomorphism for all $Y,Z \in \mathcal{P}$. We will say that $\mathcal{P}$ is \textit{locally homogeneous} if there exists a covering of $\mathcal{P}$ by open subsets $U_{a}$ of $\mathcal{B}$ such that $U_{a} \cap \mathcal{P}$ are homogeneous submanifolds of $\mathcal{B}$.
\end{definition}

As we have proved previously, the non-holonomic body-material foliation $\mathcal{NF}$ divides the body into smoothly uniform components (see theorem \ref{14.1}). We will rely on this result to provide the intuition behind the definition of homogeneity of a non-uniform body. Roughly speaking, a non-uniform Cosserat medium will be \textit{(locally) homogeneous} when each smoothly uniform material submanifold $\mathcal{NF} \left( X\right)$ is (locally) homogeneous and all the uniform material submanifolds can be ``\textit{straightened at the same time}''.
\begin{remark}
\rm
Now, suppose that $F\mathcal{B}$ is (locally) homogeneous. Then, if we take the coordinates $\left(x^{i}\right) $ on $\mathcal{B}$ given by the induced diffeomorphism $\psi$ of the definition \ref{151}, we deduce that the section $\mathcal{P}$ of Eq. (\ref{homog234}) is expressed by
\begin{equation}\label{NHOM}
\mathcal{P}\left(x^{i},y^{j}\right) = \left(\left(x^{i},y^{j}, P^{j}_{i}\right) , \delta^{j}_{i} , \dfrac{\partial P^{j}_{i} }{\partial x^{k}} + \dfrac{\partial P^{j}_{i} }{\partial y^{k}} \right)
\end{equation}
In fact, taking into account the expression of the coordinates (\ref{45}), we have that
$$ P^{j}_{i}\left(X,Y\right) \ = \ y^{j}_{i} \left( P\left(X,Y\right) \right)\left( {e_{1}}_{x}\right)$$
Hence, for a fixed two particles $X_{0},Y_{0} \in \mathcal{B}$,
\begin{eqnarray*}
P^{j}_{i,k} \left( X_{0},Y_{0} \right)  &=& \dfrac{\partial \left( y^{j}_{i} \left( P\left(X_{0},Y_{0}\right)\right) \left( {e_{1}}_{x}\right)\right)}{\partial {x^{k}} }\\
&=& \dfrac{\partial \left( y^{j}_{i} \left( P\left(X,\psi_{X_{0},Y_{0}}\left( X \right)\right)\right) \left( {e_{1}}_{x}\right)\right)}{\partial {x^{k}} }\\
&=& \dfrac{\partial P^{j}_{i} }{\partial x^{k}} + \dfrac{\partial P^{j}_{i} }{\partial y^{k}}
\end{eqnarray*}

where,
$$\psi_{X_{0},Y_{0}} = \psi^{-1}\left( \psi\left( X\right) - \psi\left( X_{0}\right)+\psi\left( Y_{0}\right)\right)$$

So, $F\mathcal{B}$ is homogeneous if we can cover $\mathcal{B}$ by local coordinate systems $\left(x^{i}\right)$ which generate (local) fields of material isomorphisms satisfying Eq. (\ref{NHOM}).
\end{remark}

Let us consider a smooth section $\mathcal{P}$ of the anchor map $\left( \overline{\alpha}, \overline{\beta}\right)$ of $\tilde{J}^{1}\left( F \mathcal{B}\right)$. Then, $\mathcal{P}$ generates a left-invariant vector field $\Theta^{\mathcal{P}}_{k}$ on $\tilde{J}^{1}\left( F \mathcal{B}\right)$ in the following way,
\begin{equation}\label{leftinvvector23}
     \Theta^{\mathcal{P}}_{k} \left( \overline{\epsilon} \left( X \right) \right) = T_{X} P^{X} \left(\dfrac{\partial }{\partial x^{k}_{\vert X}} \right)
\end{equation}
Here, for each two material points $X,Y \in \mathcal{B}$, $\mathcal{P}^{Y} \left( X \right) = \mathcal{P} \left( X,Y \right)$. If $\mathcal{P}$ is an homogeneous section, in coordinates, we have that,

\begin{equation}\label{NHOM345}
 \Theta^{\mathcal{P}}_{k} =  \dfrac{\partial }{\partial x^{k}} + \dfrac{\partial P^{j}_{i} }{\partial x^{k}} \dfrac{\partial }{\partial y^{j}_{i}} + \left[\dfrac{\partial^{2} P^{j}_{i} }{\partial x^{k}x^{l}} + \dfrac{\partial^{2} P^{j}_{i} }{\partial x^{k}y^{l}} \right] \dfrac{\partial }{\partial y^{j}_{i,l}},
\end{equation}

at the identities. So, smoothly uniform bodies are homogeneous if, and only if, the space of admissible vector fields for the couple $\left(\tilde{J}^{1} \left( F \mathcal{B} \right), \overline{\Omega} \left( \mathcal{B} \right)\right)$ may be generated by local vector fields given by Eq. (\ref{NHOM345}), for all $k$ or, in other words, we may find coordinates $\left( x^{i}\right)$ and local functions $ P^{j}_{i}$ on the macromedium $\mathcal{B}$ satisfying the equation,

\begin{equation}\label{homogeneousmaterialeqcoss2345}¡
  \dfrac{\partial  W}{\partial x^{k} } + y^{j}_{l}\dfrac{\partial P^{l}_{i} }{\partial x^{k}} \dfrac{\partial  W}{\partial y^{j}_{i} } + y^{j}_{l} \left[\dfrac{\partial^{2} P^{l}_{i} }{\partial x^{k}x^{m}} + \dfrac{\partial^{2} P^{l}_{i} }{\partial x^{k}y^{m}} \right] \dfrac{\partial  W}{\partial y^{j}_{i,m} } = 0 , \ \ k=1,2,3.
\end{equation}

Eq. (\ref{homogeneousmaterialeqcoss2345}) is, again, 
second order partial differential equation and, it will be called \textbf{homogeneity equation for Cosserat media}. Notice that, Eq. (\ref{homogeneousmaterialeqcoss2345}) makes sense even without the uniformity condition. So, we will use this to present our definition.

\begin{definition}\label{151.secondhomg}
\rm
A Cosserat medium $\mathcal{B}$ will be said to be \textit{homogeneous} if it admits coordinates $\left( x^{i}\right)$, and functions $P^{j}_{i}$, globally defined on the body, satisfying the equation,
\begin{equation}\label{2homogeneousmaterialeqcoss2345}
  \dfrac{\partial  W}{\partial x^{k} } + y^{j}_{l}\dfrac{\partial P^{l}_{i} }{\partial x^{k}} \dfrac{\partial  W}{\partial y^{j}_{i} } + \left(y^{j}_{l} \left[\dfrac{\partial^{2} P^{l}_{i} }{\partial x^{k}x^{m}} + \dfrac{\partial^{2} P^{l}_{i} }{\partial x^{k}y^{m}} \right]+ y^{j}_{l,m}\Theta^{l}_{i}\right) \dfrac{\partial  W}{\partial y^{j}_{i,m} } = 0 .
\end{equation}
for all $k = 1 , \dots , \text{dim}\left( A\overline{\Omega}\left( \mathcal{B}\right)\right)$. $\mathcal{B}$ is said to be \textit{locally homogeneous} if there exists a covering of $\mathcal{B}$ by homogeneous open sets.
\end{definition}

\begin{theorem}\label{31.cosserat}
Let $F\mathcal{B}$ be a Cosserat body. $\mathcal{B}$ will be \textit{locally homogeneous} if, and only if, for all point $X \in \mathcal{B}$ there exists a local configuration $\Psi$ of $\mathcal{B}$, defined on $F\mathcal{U}$ with $X \in \mathcal{U}$, which satisfies that
$$ j_{Y,Z}^{1} \left( \Psi^{-1} \circ \tau_{ \psi\left( Z \right) - \psi \left( Y\right)} \circ \Psi \right),$$
is a material isomorphism for all $Z \in \mathcal{U} \cap \mathcal{NF}\left( Y \right)$.
\begin{proof}
Let us assume that $F\mathcal{B}$ is (locally) homogeneous. Then, it admits coordinates $\left( x^{i}\right)$, and functions $P^{j}_{i}$, globally defined on the body, satisfying Equation (\ref{2homogeneousmaterialeqcoss2345}),
\begin{equation*}
  \dfrac{\partial  W}{\partial x^{k} } + y^{j}_{l}\dfrac{\partial P^{l}_{i} }{\partial x^{k}} \dfrac{\partial  W}{\partial y^{j}_{i} } + \left(y^{j}_{l} \left[\dfrac{\partial^{2} P^{l}_{i} }{\partial x^{k}x^{m}} + \dfrac{\partial^{2} P^{l}_{i} }{\partial x^{k}y^{m}} \right]+ y^{j}_{l,m}\Theta^{l}_{i}\right) \dfrac{\partial  W}{\partial y^{j}_{i,m} } = 0 .
\end{equation*}
for all $k = 1 , \dots , \text{dim}\left( A\overline{\Omega}\left( \mathcal{B}\right)\right)$. Thus, the local section given by
$$
\overline{\mathcal{P}}\left(x^{i},y^{j}\right) = \left(\left(x^{i},y^{j}, P^{j}_{i}\right) , \delta^{j}_{i} , \dfrac{\partial P^{j}_{i} }{\partial x^{k}} + \dfrac{\partial P^{j}_{i} }{\partial y^{k}} \right)
$$


satisfies that, its induced left-invariant vector field $\Theta^{\overline{\mathcal{P}}}_{k}$ on $\tilde{J}^{1}\left( F \mathcal{B}\right)$ is an admissible vector field for the couple $\left( \tilde{J}^{1}\left( F \mathcal{B}\right) , \overline{\Omega} \left( \mathcal{B}\right)\right)$, for all $k = 1 , \dots , \text{dim}\left( A\overline{\Omega}\left( \mathcal{B}\right)\right)$. In other words, the vector fields $\Theta^{\overline{\mathcal{P}}}_{k}$ are tangent to the non-holonomic material distribution of second order $A \overline{\Omega}\left( \mathcal{B}\right)^{T}$. Then, by Eq. (\ref{leftinvvector23}), we may assume that $\overline{\mathcal{P}}$ satisfies that
$$ \overline{\mathcal{P}} \left( Y,X \right) \in \overline{\mathcal{NF}}\left( \overline{\epsilon}\left(X\right)\right),$$
for any two material particles $X$ and $Y$, at the domain $\mathcal{U} \times \mathcal{U}$ of $\overline{\mathcal{P}}$, at the same fibre ($Y \in \mathcal{NF}\left( X \right)$), i.e., $\overline{\mathcal{P}}$ is a smooth field of material isomorphisms when it is restricted to the leaves of the non-holonomic body-material foliation of second order (recall that the non-holonomic body-material foliation $\mathcal{NF}$ divides the body $\mathcal{B}$ into maximal smoothly uniform Cosserat submanifolds).\\
Let us now fix $\overline{Z}_{0} \in \overline{F}^{2}\mathcal{B}$ such that $\overline{\pi}^{2}\left(\overline{Z}_{0}\right) = Z_{0}\in \mathcal{B}$. Then, we may define a section $\overline{\mathcal{P}}_{Z_{0}}$ of $\overline{\pi}^{2}$ given by
$$ \overline{\mathcal{P}}_{Z_{0}} \left( X \right) = \overline{\mathcal{P}} \left( Z_{0},X\right) \cdot \overline{Z}_{0} $$
Observe that, for any two particles $X,Y \in \mathcal{U}$,
\begin{equation}\label{nonhoinver}
 \overline{\mathcal{P}} \left( X , Y\right)    = \overline{\mathcal{P}}_{Z_{0}} \left( Y\right) \cdot \left[ \overline{\mathcal{P}}_{Z_{0}} \left( X \right) \right]^{-1}
\end{equation}
On the other hand, by taking into account the projections $\tilde{\pi}_{1}^{2}$ and $\overline{\pi}_{1}^{2}$, we may projects $\overline{\mathcal{P}}_{Z_{0}}$ into two sections of $F\mathcal{B}$ as follows,
\begin{itemize}
    \item $ \mathcal{Q}_{Z_{0}}  = \tilde{\pi}_{1}^{2} \circ \overline{\mathcal{P}}_{Z_{0}} $
    \item $ \mathcal{P}_{Z_{0}}  = \tilde{\pi}_{1}^{2} \circ \overline{\mathcal{P}}_{Z_{0}} $
\end{itemize}
Notice that, by the local coordinates of $\overline{\mathcal{P}}$, we have that, there are local coordinates $\psi$ of $\mathcal{B}$ in such a way that,
$$ \mathcal{Q}_{Z_{0}} \left(X\right) = j^{1}_{0,X}\left( \psi^{-1}\circ \tau_{\psi\left(X\right)}\right)$$
Thus, we construct the following map
$$
\begin{array}{rccl}
\Psi: & F\mathcal{V} & \rightarrow & F\mathcal{U} \\
& j^{1}_{0,Z}f  &\mapsto &  \mathcal{P}_{Z_{0}}\left( \psi^{-1}\left( Z\right) \right)\cdot j^{1}_{0,0}\left( \tau_{-Z}\circ f\right)
\end{array}
$$
$\Psi$ is a local principal bundle isomorphism over $\psi^{-1}$ with inverse given by
$$ j_{0,Z}^{1}g \in F\mathcal{U} \mapsto j_{0, \psi \left(Z\right)}^{1} \tau_{\psi \left(Z\right)} \cdot [\mathcal{P}_{Z_{0}} \left(Z\right)]^{-1} \cdot j_{0,Z}^{1} g.$$
In fact, we may prove that
$$ \overline{\mathcal{P}}_{Z_{0}}\left( X\right) =  j_{e_{1},\overline{X}}^{1} \left( \Psi \circ F\tau_{\psi \left(X\right)} \right)$$
Therefore, by using Eq. (\ref{nonhoinver}), we have that
$$\overline{\mathcal{P}}\left( X,Y\right) =  j_{X,Y}^{1} \left( \Psi \circ F\tau_{\psi \left(Y\right) - \psi \left(X\right)} \circ \Psi^{-1}\right), \ \forall X,Y \in \mathcal{U}.$$
Consequently, $\overline{\mathcal{P}}$ is a homogeneous sections for all the leaves of the non-holonomic body-material foliation of second order with non-empty intersection with the domain of $\overline{\mathcal{P}}$.\\
The converse is analagous.

\end{proof}
\end{theorem}

It is remarkable that, this theorem provides an intuitive view of the definition of homogeneity for non-uniform Cosserat bodies. In fact, roughly speaking, a Cosserat medium will be homogeneous, if all the leaves $\mathcal{NF} \left( X \right)$ of the unique division of the material into smoothly uniform Cosserat submanifolds are homogeneous (in the sense of Definition \ref{homogensubm}). Notice also that, the condition that all the leaves $\mathcal{NF} \left( X \right)$ are homogeneous is not enough in order to have the homogeneity of the body $\mathcal{B}$ because there is also a condition of \textit{compatibility} with the foliation structure of $\overline{\mathcal{NF}}$ given by the fact of that for all the leaves there are configuration which are homogeneous for all the leaves \textit{at the same time}. In mathematical terminology, this fact means that the homogeneous sections are induced by foliated charts for the foliation $\overline{\mathcal{NF}}$.\\

Thus, the definition of homogeneity for a smoothly uniform Cosserat medium coincides with Definition \ref{homogensubm}.\\

\section{Conclusions}
In this paper we have dealt with a model for media with microstructure, called Cosserat material. Here, we have used the so-called characteristic distributions.\\

Thus, for the case of Cosserat material, we have considered two different, but canonically defined,  characteristic distribution, called \textit{non-holonomic material distribution of second order} and \textit{holonomic material distribution of second order}, respectively. Denoting by $W$ to the mechanical response, we have proved that both distributions are characterized by the left-invariant vector fields which are in the kernel of $TW$ (Eq. (\ref{14.22})) and the complete lift of vector fields on the macromedium $\mathcal{B}$ which are in the kernel of $TW$, respectively. Therefore, we have found two different equations to construct these characteristic distributions without integrating vector fields (see definition of the characteristic distribution in section \ref{Chardistbackg94}),
\begin{itemize}
    \item \textbf{Non-holonomic material equation for Cosserat media} (\ref{materialeqcoss2345})
 \begin{equation*}
  \hspace{-1cm}  -\Theta^{i} \dfrac{\partial W }{\partial x^{i} } + \Theta^{l}_{i}\left[ y^{j}_{l}\dfrac{\partial  W}{\partial y^{j}_{i} } +  y^{j}_{l,k} \dfrac{\partial  W}{\partial y^{j}_{i,k} }\right] + \Theta^{l}_{,i} \left[ y^{j}_{,l}  \dfrac{\partial  W}{\partial y^{j}_{,i} } +  y^{j}_{m,l}\dfrac{\partial  W}{\partial y^{j}_{m,i} }\right] +   \Theta^{l}_{i,k} y^{j}_{l} \dfrac{\partial  W}{\partial y^{j}_{i,k} }= 0 
\end{equation*}

\item \textbf{Holonomic material equation for Cosserat media} (\ref{holonomicmaterialeqcoss2245})

\begin{small}
\begin{equation*}
\hspace{-1cm}-\Theta^{i} \dfrac{\partial W }{\partial x^{i} } + \dfrac{\partial  \Theta^{l}}{\partial x^{i} }\left[ y^{j}_{l}\dfrac{\partial  W}{\partial y^{j}_{i} }  +  y^{j}_{,l}  \dfrac{\partial  W}{\partial y^{j}_{,i} } + y^{j}_{l,k} \dfrac{\partial  W}{\partial y^{j}_{i,k} } +  y^{j}_{m,l}\dfrac{\partial  W}{\partial y^{j}_{m,i} }\right] +    \dfrac{\partial^{2} \Theta^{l} }{\partial x^{i}\partial x^{k} } y^{j}_{l} \dfrac{\partial  W}{\partial y^{j}_{i,k} }= 0 .
\end{equation*}
\end{small}

\end{itemize}

We have also proved that the Cosserat material is uniquely divided into smoothly uniform submanifolds and second-grade material submanifolds, respectively. Furhtermore, the space of solutions of these equations characterizes these two properties (Theorems \ref{uniformityth}, \ref{uniformityandsecondth}, \ref{uniformityandsecondth23}, and  \ref{uniformityandsecondth23456}).\\

Finally, by using these results, we have dealt we another property, homogeneity. In particular, we have been able to define, by first time, a notion of homogeneity which is valid for non-uniform materials and generalizes the known notion of homogeneity. Roughly speaking, a Cosserat material will be homogeneous, when each smoothly uniform material submanifold is homogeneous and all the uniform material submanifolds can be ``\textit{straightened at the same time}''.\\
Next, we found another differential equation,
\begin{itemize}
    \item \textbf{Homogeneity equation for Cosserat media} (\ref{2homogeneousmaterialeqcoss2345})
    \begin{equation*}
  \dfrac{\partial  W}{\partial x^{k} } + y^{j}_{l}\dfrac{\partial P^{l}_{i} }{\partial x^{k}} \dfrac{\partial  W}{\partial y^{j}_{i} } + \left(y^{j}_{l} \left[\dfrac{\partial^{2} P^{l}_{i} }{\partial x^{k}x^{m}} + \dfrac{\partial^{2} P^{l}_{i} }{\partial x^{k}y^{m}} \right]+ y^{j}_{l,m}\Theta^{l}_{i}\right) \dfrac{\partial  W}{\partial y^{j}_{i,m} } = 0,
\end{equation*}
\end{itemize}
which characterizes this intuitive notion of homogeneity.
\section*{Acknowledgments}
M. de Leon and V. M. Jiménez acknowledge the partial finantial support from MICINN Grant PID2019-106715GB-C21 and the ICMAT Severo Ochoa project CEX2019-000904-S.

\bibliographystyle{plain}

\bibliography{Library}

\begin{thebibliography}{10}

\bibitem{FGMA}
C.~M. C{\'a}mpos, M.~{Epstein}, and M.~{de Le{\'o}n}.
\newblock Functionally graded madia.
\newblock {\em International Journal of Geometric Methods in Modern Physics},
  05(03):431--455, 2008.

\bibitem{GCAPR}
G.~{Capriz}.
\newblock {\em Continua with microstructure}, volume~35 of {\em Springer Tracts
  in Natural Philosophy}.
\newblock Springer-Verlag, New York, 1989.

\bibitem{FCosse}
E.~{Cosserat} and F.~{Cosserat}.
\newblock Théorie des corps déformables.
\newblock {\em Nature}, 81(67), 1909.

\bibitem{VMMDME}
M.~{de Le\'on}, M.~{Epstein}, and V.~M. Jim\'enez.
\newblock {\em Material Geometry: Groupoids in Continuum Mechanics}.
\newblock World Scientific, Singapore, 2021.

\bibitem{SecondMat}
M.~{de Le\'on} and M.~{Esptein}.
\newblock The geometry of uniformity in second-grade elasticity.
\newblock {\em Acta Mechanica}, 114:217--224, 1996.

\bibitem{MDELAM}
M.~{de Le{\'o}n} and A.~{Mart\'\i n M\'endez}.
\newblock Principal bundle structures among second order frame bundles.
\newblock {\em Differential Geom. Appl.}, 47:202--211, 2016.

\bibitem{CELP}
C.~{Ehresmann}.
\newblock Les prolongements d'une vari\'et\'e diff\'erentiable. {V}.
  {C}ovariants diff\'erentiels et prolongements d'une structure
  infinit\'esimale.
\newblock {\em C. R. Acad. Sci. Paris}, 234:1424--1425, 1952.

\bibitem{CELC11}
C.~Ehresmann.
\newblock Introduction \`a la th\'eorie des structures infinit\'esimales et des
  pseudogroupes de {L}ie.
\newblock In {\em Colloque de topologie et g\'eom\'etrie diff\'erentielle,
  {S}trasbourg, 1952, no. 11}, page~16. La Biblioth\`eque Nationale et
  Universitaire de Strasbourg, 1953.

\bibitem{CELC12}
C.~Ehresmann.
\newblock Extension du calcul des jets aux jets non holonomes.
\newblock {\em C. R. Acad. Sci. Paris}, 239:1762--1764, 1954.

\bibitem{CELC13}
C.~Ehresmann.
\newblock Applications de la notion de jet non holonome.
\newblock {\em C. R. Acad. Sci. Paris}, 240:397--399, 1955.

\bibitem{MELZA}
M.~{El\.zanowski}, M.~{Epstein}, and J.~{\'Sniatycki}.
\newblock {$G$}-structures and material homogeneity.
\newblock {\em J. Elasticity}, 23(2-3):167--180, 1990.

\bibitem{EPST}
M.~{Epstein}.
\newblock Laminated uniformity and homogeneity.
\newblock {\em Mechanics Research Communications}, 2017.

\bibitem{MAREMDL3}
M.~{Epstein} and M.~{de Le{\'o}n}.
\newblock Homogeneity conditions for generalized {C}osserat media.
\newblock {\em J. Elasticity}, 43(3):189--201, 1996.

\bibitem{MAREMDL2}
M.~{Epstein} and M.~{de Le{\'o}n}.
\newblock Uniformity and homogeneity of elastic rods, shells and {C}osserat
  three-dimensional bodies.
\newblock {\em Arch. Math. (Brno)}, 32(4):267--280, 1996.

\bibitem{EPSTEIN1998127}
M.~Epstein and M.~{de Leon}.
\newblock Geometrical theory of uniform cosserat media.
\newblock {\em Journal of Geometry and Physics}, 26(1):127--170, 1998.

\bibitem{FGM2}
M.~{Epstein} and M.~{de Le{\'o}n}.
\newblock Homogeneity without uniformity: towards a mathematical theory of
  functionally graded materials.
\newblock {\em International Journal of Solids and Structures}, 37(51):7577 --
  7591, 2000.

\bibitem{MAREMDL4}
M.~{Epstein} and de~M.~Le{\'o}n.
\newblock The differential geometry of {C}osserat media.
\newblock 350:143--164, 1996.

\bibitem{MAREMDL}
M.~{Epstein} and de~M.~Le{\'o}n.
\newblock Geometrical theory of uniform cosserat media.
\newblock {\em J. Geom. Phys.}, 26(1-2):127--170, 1998.

\bibitem{EPSBOOK2}
M.~{Epstein} and M.~Elzanowski.
\newblock {\em Material Inhomogeneities and their Evolution: A Geometric
  Approach}.
\newblock Interaction of Mechanics and Mathematics. Springer Berlin Heidelberg,
  2007.

\bibitem{MGEOEPS}
M.~{Epstein}, V.~M. {Jim{\'e}nez}, and M.~{de Le{\'o}n}.
\newblock Material geometry.
\newblock {\em Journal of Elasticity}, 135(1):237--260, Apr 2019.

\bibitem{ERIN}
A.~C. {Eringen}.
\newblock {\em Nonlocal continuum field theories}.
\newblock Springer-Verlag, New York, 2002.

\bibitem{MD}
V.~M. {Jim{\'e}nez}, M.~{de Le{\'o}n}, and M.~{Epstein}.
\newblock Material distributions.
\newblock {\em Mathematics and Mechanics of Solids}, 25(7):1450--1458, 2017.

\bibitem{CHARDIST}
V.~M. {Jim{\'e}nez}, M.~{de Le{\'o}n}, and M.~{Epstein}.
\newblock Characteristic distribution: An application to material bodies.
\newblock {\em Journal of Geometry and Physics}, 127:19 -- 31, 2018.

\bibitem{COSVME}
V.~M. {Jim{\'e}nez}, M.~{de Le{\'o}n}, and M.~{Epstein}.
\newblock Lie groupoids and algebroids applied to the study of uniformity and
  homogeneity of cosserat media.
\newblock {\em International Journal of Geometric Methods in Modern Physics},
  15(08):1830003, 2018.

\bibitem{GENHOM}
V.~M. {Jim{\'e}nez}, M.~{de Le{\'o}n}, and M.~{Epstein}.
\newblock {\em On the homogeneity of non-uniform material bodies}, pages
  381--416.
\newblock Springer International Publishing, Berlin, 2020.

\bibitem{KONOM}
S.~I. {Kobayashi} and K.~{Nomizu}.
\newblock {\em Foundations of differential geometry. {V}ol. {I}}.
\newblock Wiley Classics Library. John Wiley \& Sons, Inc., New York, 1996.
\newblock Reprint of the 1963 original, A Wiley-Interscience Publication.

\bibitem{EKRON}
E.~{Kr{\"{o}}ner}.
\newblock {\em Mechanics of {Generalized Continua}}.
\newblock Springer, Heidelberg, 1968.

\bibitem{STRESS}
R.~Kupferman, E.~Olami, and R.~Segev.
\newblock Stress theory for classical fields.
\newblock {\em Mathematics and Mechanics of Solids}, 25(7):1472--1503, 2020.

\bibitem{KMG}
K.~C.~H. {Mackenzie}.
\newblock {\em General theory of {L}ie groupoids and {L}ie algebroids}, volume
  213 of {\em London Mathematical Society Lecture Note Series}.
\newblock Cambridge University Press, Cambridge, 2005.

\bibitem{GAMAU1}
G.~A. {Maugin}.
\newblock The method of virtual power in continuum mechanics: application to
  coupled fields.
\newblock {\em Acta Mech.}, 35(1-2):1--70, 1980.

\bibitem{GAMAU3}
G.~A. {Maugin}.
\newblock On the structure of the theory of polar elasticity.
\newblock {\em R. Soc. Lond. Philos. Trans. Ser. A Math. Phys. Eng. Sci.},
  356(1741):1367--1395, 1998.

\bibitem{WNOLLTHE}
W.~Noll.
\newblock {\em On the continuity of the solid and fluid states}.
\newblock ProQuest LLC, Ann Arbor, MI, 1954.
\newblock Thesis (Ph.D.)--Indiana University.

\bibitem{WNOLL}
W.~{Noll}.
\newblock Materially uniform simple bodies with inhomogeneities.
\newblock {\em Arch. Rational Mech. Anal.}, 27:1--32, 1967.

\bibitem{SAUND}
D.~J. Saunders.
\newblock {\em The Geometry of Jet Bundles}.
\newblock London Mathematical Society Lecture Note Series. Cambridge University
  Press, 1989.

\bibitem{PS}
P.~{Stefan}.
\newblock Accessible sets, orbits, and foliations with singularities.
\newblock {\em Proc. London Math. Soc. (3)}, 29:699--713, 1974.

\bibitem{HJS}
H.~J. {Sussmann}.
\newblock Orbits of families of vector fields and integrability of
  distributions.
\newblock {\em Trans. Amer. Math. Soc.}, 180:171--188, 1973.

\bibitem{CTRUE}
C.~{Truesdell} and W.~{Noll}.
\newblock {\em The non-linear field theories of mechanics}.
\newblock Springer-Verlag, Berlin, third edition, 2004.
\newblock Edited and with a preface by Stuart S. Antman.

\end{thebibliography}

\end{document}